\documentclass[11pt]{amsart}
\usepackage{a4wide,amsmath,amssymb,amscd,verbatim,latexsym}

\newtheorem{prop}[equation]{Proposition}

\newtheorem{thm}[equation]{Theorem}
\newtheorem{cor}[equation]{Corollary}
\newtheorem{lem}[equation]{Lemma}
\newtheorem{exa}[equation]{Example}
\newtheorem{exas}[equation]{Examples}

\theoremstyle{definition}

\newtheorem{defn}[equation]{Definition}
\newtheorem{defns}[equation]{Definitions}
\newtheorem{cons}[equation]{Construction}

\numberwithin{equation}{section}

\newcommand{\spandsp}{\mbox{$\qquad\text{and}\qquad$}}

\newcommand{\sands}{\mbox{$\quad\text{and}\quad$}}

\newcommand{\Hom}{\operatorname{Hom}}
\newcommand{\Tor}{\operatorname{Tor}}

\newcommand{\col}{\operatorname{colim}}
\newcommand{\coltmg}{\operatorname{colim^{\scat{tmg}}}}
\newcommand{\coltmon}{\operatorname{colim^{\scat{tmon}}}}

\newcommand{\hoc}{\operatorname{hocolim}}
\newcommand{\hoctmg}{\operatorname{hocolim^{\scat{tmg}}}}

\newcommand{\vble}{\;\,}

\newcommand{\bC}{\mathbb{C}}

\newcommand{\bE}{\mathbb{E}}
\newcommand{\bF}{\mathbb{F}}

\newcommand{\bR}{\mathbb{R}}
\newcommand{\bZ}{\mathbb{Z}}

\newcommand{\bbox}{\mathbin{\Box}}

\newcommand{\cat}[1]{\mbox{\sc #1}}
\newcommand{\scat}[1]{\mbox{\scriptsize{\sc #1}}}
\newcommand{\sdel}{\mbox{\scriptsize{$\Delta$}}}
\newcommand{\del}{\mbox{\footnotesize{$\Delta$}}}

\newcommand{\lca}{\begin{picture}(40,7.5)\put(6,-2){$\relbar\joinrel\hspace{-1pt}\relbar\joinrel\hspace{-1pt}\longrightarrow$}\put(6,4){$\relbar\joinrel\hspace{-1pt}\relbar\joinrel\hspace{-1pt}\longrightarrow$}
\end{picture}}
\newcommand{\coeq}{\mathop{\lca}}
\newcommand{\Sin}{\mbox{\it Sin\/}}
\newcommand{\fcat}[2]{\mbox{\raisebox{1pt}{\rm\scriptsize [}{\sc #1},{\sc #2}\raisebox{1pt}{\rm\scriptsize ]}}\hspace{1pt}}

\newcommand{\srf}[1]{\mbox{${\text{\it SR}^F}$}}
\newcommand{\under}{\!\downarrow\!}

\newcommand{\SR}{\mbox{\it SR\/}}
\newcommand{\daja}{Davis-Januszkiewicz}
\newcommand{\daaja}{Davis and Januszkiewicz}
\newcommand{\djs}{\mbox{\it DJ\/}}
\newcommand{\Fl}{\mbox{\it Fl\/}}
\newcommand{\Con}{\mbox{\it Con\/}}

\newcommand{\DP}{\mbox{\it DP\/}}
\newcommand{\Art}{\mbox{\it Art\/}}
\newcommand{\Cox}{\mbox{\it Cox\/}}
\newcommand{\Cir}{\mbox{\it Cir\/}}

\newcommand{\bigfree}{\mathop{\mbox{\Large$*$}}}

\begin{document}
\bibliographystyle{plain}
\title{Colimits, Stanley-Reisner algebras, and loop spaces}
\author{Taras Panov}
\address{Department of Mathematics and Mechanics, Moscow State
University, 119899 Moscow, Russia\newline
\emph{and} Institute for Theoretical and Experimental Physics,
117259 Moscow, Russia}
\email{tpanov@mech.math.msu.su}
\author{Nigel Ray}
\address{Department of Mathematics, University of Manchester,
Manchester M13~9PL, England}
\email{nige@ma.man.ac.uk}
\author{Rainer Vogt}
\address{Fachbereich Mathematik/Informatik, Universit\"{a}t Osnabr\"{u}ck,
D-49069 Osnabr\"{u}ck, Germany}
\email{rainer@mathematik.uni-osnabrueck.de}

\thanks{The first author was supported by a Royal Society/NATO
Postdoctoral Fellowship at the University of Manchester, and also by the
Russian Foundation for Basic Research, grant number 01-01-00546}


\begin{abstract}

We study diagrams associated with a finite simplicial complex $K$,
in various algebraic and topological categories. We relate their
colimits to familiar structures in algebra, combinatorics,
geometry and topology. These include: right-angled Artin and
Coxeter groups (and their complex analogues, which we call {\it
circulation groups}); Stanley-Reisner algebras and coalgebras;
Davis and Januszkiewicz's spaces {\it DJ}$(K)$ associated with
toric manifolds and their generalisations; and coordinate subspace
arrangements. When $K$ is a flag complex, we extend well-known
results on Artin and Coxeter groups by confirming that the
relevant circulation group is homotopy equivalent to the space of
loops $\varOmega${\it DJ}$(K)$. We define homotopy colimits for
diagrams of topological monoids and topological groups, and show
they commute with the formation of classifying spaces in a
suitably generalised sense. We deduce that the homotopy colimit of
the appropriate diagram of topological groups is a model for
$\varOmega${\it DJ}$(K)$ for an arbitrary complex $K$, and that
the natural projection onto the original colimit is a homotopy
equivalence when $K$ is flag. In this case, the two models are
compatible.

\end{abstract}

\maketitle

\section{Introduction}\label{intro}

In this work we study diagrams associated with a finite simplicial
complex $K$, in various algebraic and topological categories. We are
particularly interested in colimits and homotopy colimits of such
diagrams.

We are motivated by \daaja's investigation \cite{daja:cpc} of toric
manifolds, in which $K$ first arises as the boundary of the quotient
polytope. In the course of their cohomological computations, \daaja\
construct real and complex versions of a space whose cohomology ring
is isomorphic to the Stanley-Reisner algebra of $K$, over $\bZ/2$ and
$\bZ$ respectively. We denote the homotopy type of these spaces by
$\djs(K)$, and follow Buchstaber and Panov \cite{bupa:tact} by
describing them as colimits of diagrams of classifying spaces. In this
context, an exterior version arises naturally as an alternative.
Suggestively, the cohomology algebras and homology coalgebras of the
$\djs(K)$ may be expressed as the limits and colimits of analogous
diagrams in the corresponding algebraic category.

When colimits of similar diagrams are taken in a category of discrete
groups, they yield right-angled Coxeter and Artin groups. These are
more usually described by a complementary construction involving only
the $1$-skeleton $K^{(1)}$ of $K$. Whenever $K$ is determined entirely
by $K^{(1)}$ it is known as a flag complex, and results such as those
of \cite{daja:cpc} and \cite{kiro:hca} may be interpreted as showing
that the associated Coxeter and Artin groups are homotopy equivalent
to the loop spaces $\varOmega\djs(K)$, in the real and exterior cases
respectively. In other words, the groups are discrete models for the
loop spaces. These observations raise the possibility of modelling
$\varOmega\djs(K)$ in the complex case, and for arbitrary $K$, by
colimits of diagrams in a suitably defined category of topological
monoids. Our primary aim is to carry out this programme.

Before we begin, we must therefore confirm that our categories are
sufficiently cocomplete for the proposed colimits to exist. We show
that this is indeed the case (as predicted by folklore), and explain
how the complex version of $\varOmega\djs(K)$ is modelled by the
colimit of a diagram of tori whenever $K$ is flag. We refer to the
colimit as a {\it circulation group\/}, and consider it as the complex
analogue of the corresponding right-angled Coxeter and Artin
groups. Of course, it is also determined by $K^{(1)}$. On the other
hand, there are simple examples of non-flag complexes for which the
colimit groups cannot possibly model $\varOmega\djs(K)$ in any of the
real, exterior, or complex cases. More subtle constructions are
required.

Since we are engaged with homotopy theoretic properties of colimits,
it is no great surprise that the appropriate model for arbitrary
complexes $K$ is a homotopy colimit. Considerable care has to be
taken in formulating the construction for topological monoids, but
the outcome clarifies the status of the original colimits when $K$
is flag; flag complexes are precisely those for which the colimit
and the homotopy colimit coincide. Our main result is therefore that
$\varOmega\djs(K)$ is modelled by the homotopy colimit of the
relevant diagram of topological groups, in all three cases and for
arbitrary $K$. When $K$ is flag, the natural projection onto the
original colimit is a homotopy equivalence, and is compatible with
the two model maps. Our proof revolves around the fact that homotopy
colimits commute with the classifying space functor, in a context
which is considerably more general than is needed here.

For particular complexes $K$, our constructions have interesting
implications for traditional homotopy theoretic invariants such as
Whitehead products, Samelson products, and their higher analogues and
iterates. We hope to deal with these issues in subsequent work
\cite{para:hhc}.

We now summarise the contents of each section.

It is particularly convenient to use the language of enriched category
theory, so we devote Section \ref{capr} to establishing the notation,
conventions and results that we need. These include a brief discussion
of simplicial objects and their realisations, and verification of the
cocompleteness of our category of topological monoids in the enriched
setting. Readers who are familiar with this material, or willing to
refer back to Section \ref{capr} as necessary, may proceed directly to
Section \ref{baco}, where we introduce the relevant categories and
diagrams associated with a simplicial complex $K$. They include
algebraic and topological examples, amongst which are the exponential
diagrams $G^K$; here $G$ denotes one of the cyclic groups $C_2$ and
$C$, or the circle group $T$, in the real, exterior, and complex cases
respectively.

We devote Section \ref{co} to describing the limits and colimits of
these diagrams. Some are identified with standard constructions such
as the Stanley-Reisner algebra of $K$ and the \daja\ spaces $\djs(K)$,
whereas the $G^K$ yield right-angled Coxeter and Artin groups, or
circulation groups respectively. In Section \ref{fihoco} we study
aspects of the diagrams involving associated fibrations and homotopy
colimits. We note connections with coordinate subspace arrangements.

We introduce the model map $f_K\colon\coltmg
G^K\rightarrow\varOmega\djs(K)$ in Section \ref{flcoco}, and determine
the connectivity of its homotopy fibre in terms of combinatorial
properties of $K$. The results confirm that $f_K$ is a homotopy
equivalence whenever $K$ is flag, and quantify its failure for general
$K$. In our final Section \ref{hocotomo} we consider suitably
well-behaved diagrams $D$ of topological monoids, and prove that the
homotopy colimit of the induced diagram of classifying spaces is
homotopy equivalent to the classifying space of the homotopy colimit
of $D$, taken in the category of topological monoids. By application
to the exponential diagrams $G^K$, we deduce that our generalised
model map $h_K\colon\hoctmg G^K\rightarrow\varOmega\djs(K)$ is a
homotopy equivalence for all complexes $K$. We note that the two models
are compatible, and homotopy equivalent, when $K$ is flag.

We take the category $\cat{top}$ of $k$-spaces $X$ and continuous
functions $f\colon X\rightarrow Y$ as our underlying topological
framework, following \cite{vo:cct}. Every function space $Y^X$ is
endowed with the corresponding $k$-topology. Many of the spaces we
consider have a distinguished basepoint $*$, and we write $\cat{top}_+$
for the category of pairs $(X,*)$ and basepoint preserving maps; the
forgetful functor $\cat{top}_+\rightarrow\cat{top}$ is faithful. For any
object $X$ of $\cat{top}$, we may add a disjoint basepoint to obtain a
based space $X_+$. The $k$-function space $(Y,*)^{(X,*)}$ has the
trivial map $X\rightarrow *$ as basepoint. In some circumstances we need
$(X,*)$ to be {\it well-pointed}, in the sense that the inclusion of the
basepoint is a closed cofibration, and we emphasise this requirement as
it arises.

Several other important categories are related to $\cat{top}_+$.
These include $\cat{tmonh}$, consisting of associative topological
monoids and homotopy homomorphisms \cite{br:shc} (essentially
equivalent to Sugawara's strongly homotopy multiplicative maps
\cite{su:hcg}), and its subcategory \cat{tmon}, in which the
homorphisms are strict. Again, the forgetful functor
$\cat{tmon}\rightarrow\cat{top}_+$ is faithful. Limiting the objects
to topological groups defines a further subcategory \cat{tgrp},
which is full in \cat{tmon}. In all three cases the identity element
$e$ is the basepoint, and we may sometimes have to insist that
objects are well-pointed. The Moore loop space $\varOmega X$ is a
typical object in \cat{tmonh} for any pair $(X,*)$, and the
canonical inclusion $M\rightarrow\varOmega BM$ is a homotopy
homomorphism for any well-pointed topological monoid $M$.

For each $m\geq 0$ we consider the small categories $\cat{id}(m)$,
which consist of $m$ objects and their identity morphisms; in
particular, we use the based versions $\cat{id}_\varnothing(m)$,
which result from adjoining an initial object $\varnothing$. Given a
topological monoid $M$, the associated topological category
$\cat{c}(M)$ consists of one object, and one morphism for each
element of $M$. Segal's \cite{se:css} classifying space
$B\cat{c}(M)$ then coincides with the standard classifying space
$BM$.

Given objects $X_0$ and $X_n$ of any category $\cat{c}$, we
denote the set of $n$-composable morphisms
\[
X_0\stackrel{f_1}{\longrightarrow}X_1\stackrel{f_2}{\longrightarrow}
\dots
\stackrel{f_n}{\longrightarrow}X_n
\]
by $\cat{c}_n(X_0,X_n)$, for all $n\geq 0$. Thus $\cat{c}_1(X,Y)$ is
the morphism set $\cat{c}(X,Y)$ for all objects $X$ and $Y$, and
$\cat{c}_0(X,X)$ consists solely of the identity morphism on $X$.

In order to distinguish between them, we write $T$ for the
multiplicative topological group of unimodular complex numbers, and
$S^1$ for the circle. Similarly, we discriminate between the cyclic
group $C_2$ and the ring of residue classes $\bZ/2$, and between the
infinite cyclic group $C$ and the ring of integers $\bZ$.

The first and second authors benefitted greatly from illuminating
discussions with Bill Dwyer at the International Conference on Algebraic
Topology held on the Island of Skye during June 2001. They are
particularly grateful to the organisers for providing the opportunity to
work in such magnificent surroundings.

\section{Categorical Prerequisites}\label{capr}

We refer to the books of Kelly \cite{ke:bce} and Borceux \cite{bo:hca}
for notation and terminology associated with the theory of enriched
categories, and to Barr and Wells \cite{bawe:ttt} for background on
the theory of monads (otherwise known as triples). For more specific
results, we cite \cite{elkrmama:rma} and \cite{hovo:mts}. Unless
otherwise stated, we assume that all our categories are enriched in
one of the topological senses described below, and that functors are
continuous. In many cases the morphism sets are finite, and therefore
invested with the discrete topology.

Given an arbitrary category $\cat{r}$, we refer to a covariant functor
$D\colon\cat{a}\rightarrow\cat{r}$ as an {\it $\cat{a}$-diagram} in
$\cat{r}$, for any small category $\cat{a}$. Such diagrams are the
objects of a category $\fcat{a}{r}$, whose morphisms are natural
transformations of functors. We may interpret any object $X$ of
$\cat{r}$ as a constant diagram, which maps each object of $\cat{a}$
to $X$ and every morphism to the identity.
\begin{exas}\label{simco}
Let $\del$ be the category whose objects are the ordinals
$(n)=\{0,1,\dots,n\}$, where $n\geq 0$, and whose morphisms are the
nondecreasing functions; then $\del^{op}$- and $\del$-diagrams are
simplicial and cosimplicial objects of $\cat{a}$ respectively. In
particular, $\varDelta\colon\del\rightarrow\cat{top}$ is the
cosimplicial space which assigns the standard $n$-simplex
$\varDelta(n)$ to each object $(n)$. Its pointed analogue
$\varDelta_+$ is given by $\varDelta_+(n)=\varDelta(n)_+$.

If $M$ is a topological monoid, then $\cat{c}(M)$- and
$\cat{c}(M)^{op}$-diagrams in $\cat{top}$ are left and right $M$-spaces
respectively.
\end{exas}

We recall that $(\cat{s},\Box,\varPhi)$ is a {\it symmetric
monoidal\/} category if the bifunctor
$\bbox\colon\cat{s}\times\cat{s}\rightarrow\cat{s}$ is coherently
associative and commutative, and $\varPhi$ is a coherent unit object.
Such an $\cat{s}$ is {\it closed\/} if there is a bifunctor
$\cat{s}\times\cat{s}^{op}\rightarrow\cat{s}$, denoted by
$(Z,Y)\mapsto[Y,Z]$, which satisfies the adjunction
\[
\cat{s}(X\bbox Y,Z)\cong\cat{s}(X,[Y,Z])
\]
for all objects $X$, $Y$, and $Z$ of $\cat{s}$. A category $\cat{r}$ is
{\it $\cat{s}$-enriched\/} when its morphism sets are identified with
objects of $\cat{s}$, and composition factors naturally through
$\Box$. A closed symmetric monoidal category is canonically
self-enriched, by identifying $\cat{s}(X,Y)$ with $[X,Y]$. Henceforth,
$\cat{s}$ denotes such a category.

\begin{exa}\label{rrdiag}
Any small $\cat{s}$-enriched category $\cat{a}$ determines a diagram
$A\colon\cat{a}\times\cat{a}^{op}\rightarrow\cat{s}$, whose value at
$(a,b)$ is the morphism object $\cat{a}(b,a)$.
\end{exa}

An {\it $\cat{s}$-functor} $\cat{q}\rightarrow\cat{r}$ of
$\cat{s}$-enriched categories acts on morphism sets as a morphism of
$\cat{s}$. The category $\fcat{q}{r}$ of such functors has morphisms
consisting of natural transformations, and is also
$\cat{s}$-enriched. The $\cat{s}$-functors
$F\colon\cat{q}\rightarrow\cat{r}$ and
$U\colon\cat{r}\rightarrow\cat{q}$ are {\it $\cat{s}$-adjoint} if
there is a natural isomorphism
\[
\cat{r}(F(X),Y)\cong\cat{q}(X,U(Y))
\]
in $\cat{s}$, for all objects $X$ of $\cat{q}$ and $Y$ of $\cat{r}$.
\begin{exas}
The categories $\cat{top}$ and $\cat{top}_+$ are symmetric monoidal
under cartesian product $\times$ and smash product $\wedge$
respectively, with unit objects the one-point space $*$ and the
zero-sphere $*_+$. Both are closed, and therefore self-enriched, by
identifying $[X,Y]$ with $Y^X$ and $(Y,*)^{(X,*)}$ respectively.
\end{exas}
Since $(Y,*)^{(X,*)}$ inherits the subspace topology from $Y^X$, the
induced $\cat{top}$-enrichment of $\cat{top}_+$ is compatible with its
self-enrichment. Both $\cat{tmon}$ and $\cat{tgrp}$ are
$\cat{top}_+$-enriched by restriction.

In certain situations it is helpful to reserve the notation $\cat{t}$
for either or both of the self-enriched categories $\cat{top}$ and
$\cat{top}_+$. Similarly, we reserve $\cat{tmg}$ for either or both of
the $\cat{top}_+$-enriched categories $\cat{tmon}$ and $\cat{tgrp}$.

It is well known that $\cat{top}$ and $\cat{top}_+$ are complete and
cocomplete, in the standard sense that every small diagram has a limit
and colimit. Completeness is equivalent to the existence of products
and equalizers, and cocompleteness to the existence of coproducts and
coequalizers. Both $\cat{top}$ and $\cat{top}_+$ actually admit {\it
indexed limits\/} and {\it indexed colimits} \cite{ke:bce}, involving
topologically parametrized diagrams in the enriched setting; in other
words, $\cat{t}$ is {\it $\cat{t}$-complete\/} and {\it
$\cat{t}$-cocomplete}. A summary of the details for $\cat{top}$ can be
found in \cite{mcscvo:thh}.

Amongst indexed limits and colimits, the enriched analogues of products
and coproducts are particularly important.
\begin{defns}\label{tenscotens}
An $\cat{s}$-enriched category $\cat{r}$ is {\it tensored\/}
and {\it cotensored\/} over $\cat{s}$ if there exist bifunctors
$\cat{r}\times\cat{s}\rightarrow\cat{r}$ and
$\cat{r}\times\cat{s}^{op}\rightarrow\cat{r}$ respectively, denoted by
\[
(X,Y)\longmapsto X\otimes Y\sands(X,Y)\longmapsto X^Y,
\]
together with natural isomorphisms
\begin{equation}\label{tcotads}
\cat{r}(X\otimes
Y,Z)\cong\cat{s}(Y,\cat{r}(X,Z))\cong\cat{r}(X,Z^Y)
\end{equation}
in $\cat{s}$, for all objects $X$, $Z$ of $\cat{r}$ and $Y$ of
$\cat{s}$.
\qed
\end{defns}
For any such $\cat{r}$, there are therefore natural isomorphisms
\begin{equation}\label{nats}
X\otimes\varPhi\cong X\cong X^\varPhi\sands
X\otimes(Y\bbox W)\cong(X\otimes Y)\otimes W.
\end{equation}
Every $\cat{s}$ is tensored over itself by $\Box$, and cotensored by
$[\vble,\vble]$.
\begin{exas}\label{toptenco}
The categories $\cat{t}$ are tensored and cotensored over themselves;
so $X\otimes Y$ and $X^Y$ are given by $X\times Y$ and $X^Y$ in
$\cat{top}$, and by $X\wedge Y$ and $(Y,*)^{(X,*)}$ in $\cat{top}_+$.
\end{exas}

The r\^ole of tensors and cotensors is clarified by the following
results of Kelly \cite[(3.69)-(3.73)]{ke:bce}. Here and henceforth, we
take $\cat{s}$ to be complete and cocomplete in the standard sense.
\begin{thm}\label{kelly}
An $\cat{s}$-enriched category is $\cat{s}$-complete if and only if it
is complete, and cotensored over $\cat{s}$; it is $\cat{s}$-cocomplete
if and only if it is cocomplete, and tensored over $\cat{s}$.
\end{thm}
Theorem \ref{kelly} asserts that standard limits and colimits may
themselves be enriched in the presence of tensors and cotensors, since
they are special cases of indexed limits and colimits. Given an
$\cat{a}$-diagram $D$ in $\cat{r}$, where $\cat{a}$ is also
$\cat{s}$-enriched, we deduce that the natural bijections
\begin{equation}\label{eneqcoeq}
\cat{r}(X,\lim D)\longleftrightarrow\fcat{a}{r}(X,D)\spandsp
\cat{r}(\col D,Y)\longleftrightarrow\fcat{a}{r}(D,Y)
\end{equation}
are isomorphisms in $\cat{s}$, for any objects $X$ and $Y$ of $\cat{r}$.

It is convenient to formulate several properties of $\cat{tmon}$ and
$\cat{tgrp}$ by observing that both categories are
$\cat{top}_+$-complete and -cocomplete. We appeal to the monad
associated with the forgetful functor
$U\colon\cat{tmg}\rightarrow\cat{top}_+$; in both cases it has a left
$\cat{top}_+$-adjoint, given by the free monoid or free group functor
$F$. The composition $U\cdot F$ defines a $\cat{top}_+$-monad
$L\colon\cat{top}_+\rightarrow\cat{top}_+$, whose category
$\cat{top}_+^L$ of algebras is precisely $\cat{tmg}$. We write $V$ for
the forgetful functor $\cat{tgrp}\rightarrow\cat{tmon}$, whose left
$\cat{top}_+$-adjoint is the universal group functor.
\begin{prop}\label{tmontcc}
The categories $\cat{tmon}$ and $\cat{tgrp}$ are $\cat{top}_+$-complete
and -cocomplete; moreover, $V$ preserves indexed colimits.
\end{prop}
\begin{proof}
We consider the forgetful functor
$\cat{top}_+^L\rightarrow\cat{top}_+$, noting that
$\cat{top}_+$ is $\cat{top}_+$-complete by Theorem \ref{kelly}.

Part (i) of \cite[VII, Proposition 2.10]{elkrmama:rma} asserts that
the forgetful functor creates all indexed limits, confirming that
$\cat{tmg}$ is $\cat{top}_+$-complete. Part (ii) asserts that
$\cat{top}_+^L$ is $\cat{top}_+$-cocomplete if $L$ preserves reflexive
coequalizers, which need only be verified for $U$ because $F$
preserves colimits. The result follows for an arbitrary reflexive pair
$(f,g)$ in $\cat{tmg}$ by using the right inverse to show that the
coequalizer of $(U(f),U(g))$ in $\cat{top}_+$ is itself in the image
of $U$, and lifts to the coequalizer of $(f,g)$.

Finally, $V$ preserves indexed colimits because it is left adjoint to
a $\cat{top}_+$-functor $\cat{tmon}\rightarrow\cat{tgrp}$, which
associates to each topological monoid its subgroup of invertible
elements.
\end{proof}

In view of Proposition \ref{tmontcc} we may form colimits of diagrams
in $\cat{tmg}$ by applying $\coltmon$, even when the diagram consists
entirely of topological groups. Pioneering results on the completeness
and cocompleteness of categories of topological monoids and
topological groups may be found in \cite{brha:tgo}.

Our main deduction from Proposition \ref{tmontcc} is that $\cat{tmon}$
and $\cat{tgrp}$ are tensored over $\cat{top}_+$. By studying the
isomorphisms \eqref{tcotads}, we may construct the tensors explicitly;
they are described as pushouts in \cite[2.2]{scvo:cae}.
\begin{cons}\label{tmontens}
For any objects $M$ of $\cat{tmon}$ and $Y$ of $\cat{top}_+$, the
{\it tensored monoid} $M\circledast Y$ is the quotient of the free
topological monoid on $U(M)\wedge Y$ by the relations
\[
(m,y)(m',y)=(mm',y)\quad\text{for all $m,\,m'\in M$ and $y\in Y$}.
\]
For any object $G$ of $\cat{tgrp}$, the {\it tensored group}
$G\circledast Y$ is the topological group $V(G)\circledast Y$.

The {\it cotensored monoid} $M^Y$ and {\it cotensored group} $G^Y$ are
the function spaces $\cat{top}_+(Y,M)$ and $\cat{top}_+(Y,G)$
respectively, under pointwise multiplication.
\qed
\end{cons}

Given a category $\cat{r}$ which is tensored and cotensored over
$\cat{s}$, we may now describe several categorical constructions. They
are straightforward variations on \cite[2.3]{hovo:mts}, and initially
involve three diagrams. The first
is $D\colon\cat{b}^{op}\rightarrow\cat{r}$, the second
$E\colon\cat{b}\rightarrow\cat{s}$, and the third
$F\colon\cat{b}\rightarrow\cat{r}$.
\begin{defns}\label{tprodfus}
The {\it tensor product} $D\otimes_{\scat{b}}E$ is the coequalizer of
\[
\coprod_{g:b_0\to b_1}\hspace{-2mm}D(b_1)\otimes E(b_0)
\coeq\limits_\beta^\alpha \coprod_b D(b)\otimes E(b)
\]
in $\cat{r}$, where $g$ ranges over the morphisms of $\cat{b}$, and
$\alpha|_{g}=D(g)\otimes 1$ and $\beta|_{g}=1\otimes E(g)$. The {\it
homset} $\Hom_{\scat{b}}(E,F)$ is the equalizer of
\[
\prod_bF(b)^{E(b)} \coeq\limits_\beta^\alpha \prod_{g:b_0\to
b_1}\hspace{-2mm} F(b_1)^{E(b_0)}
\]
in $\cat{r}$, where $\alpha=\prod_g\cdot E(g)$ and
$\beta=\prod_g F(g)\cdot$\;.
\qed
\end{defns}
We may interpret the elements of $\Hom_{\scat{b}}(E,F)$ as mappings
from the diagram $E$ to the diagram $F$, using the cotensor pairing.

\begin{exas}\label{tencotexa}
Consider the case $\cat{r}=\cat{s}=\cat{top}$ or $\cat{top}_+$, with
$\cat{b}=\del$. Given simplicial spaces
$X_\bullet\colon\del^{op}\rightarrow\cat{top}$ and
$Y_\bullet\colon\del^{op}\rightarrow\cat{top}_+$, the tensor products
\[
|X_\bullet|=X_\bullet\times_{\sdel}\varDelta
\spandsp
|Y_\bullet|=Y_\bullet\wedge_{\sdel}\varDelta_+
\]
represent their topological realisation \cite{ma:soa} in $\cat{top}$
and $\cat{top}_+$ respectively. If we choose $\cat{r}=\cat{tmg}$ and
$\cat{s}=\cat{top}_+$, a simplicial object
$M_\bullet\colon\del^{op}\rightarrow\cat{tmg}$ has internal and
topological realisations
\[
|M_\bullet|_{\scat{tmg}}=M_\bullet\circledast_{\sdel}\varDelta_+
\spandsp
|M_\bullet|=U(M)_\bullet\wedge_{\sdel}\varDelta_+
\]
in $\cat{tmg}$ and $\cat{top}_+$ respectively. Since $|\vble|$
preserves products, $|M_\bullet|$ actually lies in $\cat{tmg}$.

If $\cat{r}=\cat{s}$, then $D\otimes_{\scat{b}}\varPhi$ is $\col D$,
where $\varPhi$ is the trivial $\cat{b}$-diagram. Also,
$\Hom_{\scat{b}}(E,F)$ is the morphism set $\fcat{b}{r}(E,F)$,
consisting of the natural transformations $E\rightarrow F$.
\end{exas}
For $Y_\bullet$ in Examples \ref{tencotexa}, its $\cat{top}$- and
$\cat{top}_+$-realisations are homeomorphic because basepoints of the
$Y_n$ represent degenerate simplices for $n>0$. We identify
$|M_\bullet|_{\scat{tmg}}$ with $|M_\bullet|$ in Section
\ref{hocotomo}.

We need certain generalisations of Definitions \ref{tprodfus}, in
which analogies with homological algebra become apparent. We extend
the first and second diagrams to
$D\colon\cat{a}\times\cat{b}^{op}\rightarrow\cat{r}$ and
$E\colon\cat{b}\times\cat{c}^{op}\rightarrow\cat{s}$, and replace the
third by
$F\colon\cat{c}\times\cat{d}^{op}\rightarrow\cat{s}$ or
$G\colon\cat{a}\times\cat{c}^{op}\rightarrow\cat{r}$. Then
$D\otimes_{\scat{b}}E$ becomes an
$(\cat{a}\times\cat{c}^{op})$-diagram in $\cat{r}$, and
$\Hom_{\scat{c}^{op}}(E,G)$ becomes an
$(\cat{a}\times\cat{b}^{op})$-diagram in $\cat{r}$. The extended
diagrams reduce to the originals by judicious substitution, such as
$\cat{a}=\cat{c}=\cat{id}$ in $D$ and $E$.

\begin{exa}\label{tosincom}
Consider the case $\cat{r}=\cat{s}=\cat{top}_+$, with
$\cat{a}=\cat{c}=\cat{id}$ and $\cat{b}=\del$. Given $E=\varDelta_+$
as before, and $G$ a constant diagram
$Z\colon\cat{id}\rightarrow\cat{top}_+$, then
$\Hom_{\scat{c}^{op}}(E,G)$ coincides with the total singular complex
$\Sin(Z)$ as an object of $\fcat{\del$^{op}$}{top$_+$}$. If
$\cat{r}=\cat{tmg}$ and $N\colon\cat{id}\rightarrow\cat{tmg}$ is a
constant diagram, then $\Sin(N)$ is an object of
$\fcat{\del$^{op}$}{tmg}$.
\end{exa}

Important properties of tensor products are described by the natural
equivalences
\begin{equation}\label{asscot}
D\otimes_{\scat{b}}B\cong D\spandsp
(D\otimes_{\scat{b}}E)\otimes_{\scat{c}}F\cong
D\otimes_{\scat{b}}(E\otimes_{\scat{c}}F)
\end{equation}
of $(\cat{a}\times\cat{b}^{op})$- and
$(\cat{a}\times\cat{d}^{op})$-diagrams respectively, in $\cat{r}$. The
first equivalence applies Example \ref{rrdiag} with $\cat{a}=\cat{b}$,
and the second uses the isomorphism of \eqref{nats}. The adjoint
relationship between $\otimes$ and $\Hom$ is expressed by the equivalences
\begin{equation}\label{duht}
\fcat{a$\times$c$^{op}$}{r}(D\otimes_{\scat{b}}E,G)\cong
\fcat{b$\times$c$^{op}$}{s}(E,\fcat{a}{r}(D,G))\cong
\fcat{a$\times$b$^{op}$}{r}(D,\Hom_{\scat{c}^{op}}(E,G)),
\end{equation}
which extend the tensor-cotensor relations \eqref{tcotads}, and are a
consequence of the constructions.

\begin{exas}\label{retoad}
Consider the data of Example {\rm \ref{tosincom}}, and suppose that
$D$ is a simplicial pointed space
$Y_\bullet\colon\del^{op}\rightarrow\cat{top}_+$. Then the adjoint
relation \eqref{duht} provides a homeomorphism
\[
\cat{top}_+(|Y_\bullet|,Z)\cong
\fcat{\del$^{op}$}{top$_+$}(Y_\bullet,\Sin(Z)).
\]
If $\cat{r}=\cat{tmg}$ and $\cat{s}=\cat{top}_+$, and $M_\bullet$ is a
simplicial object in $\cat{tmg}$, we obtain a homeomorphism
\[
\cat{tmg}(|M_\bullet|_{\scat{tmg}},N)\cong
\fcat{\del$^{op}$}{tmg}(M_\bullet,\Sin(N))
\]
for any object $N$ of $\cat{tmg}$.

If $\cat{r}=\cat{s}$ and $E=\varPhi$, the relations \eqref{duht}
reduce to the second isomorphism \eqref{eneqcoeq}.
\end{exas}
The first two examples extend the classic adjoint relationship between
$|\vble|$ and $\Sin$.

We now assume $\cat{r}=\cat{s}=\cat{top}$. We let $D$ be an
$(\cat{a}\times\cat{b}^{op})$-diagram as above, and define
$B_\bullet(*,\cat{a},D)$ to be a degenerate form of the {\it 2-sided
bar construction}. It is a $\cat{b}^{op}$-diagram of simplicial
spaces, given as a
$\cat{b}^{op}\times\del^{op}$-diagram in $\cat{top}$ by
\begin{equation}\label{tsbar}
(b,(n))\longmapsto
\bigsqcup_{a_0,a_n}D(b,a_0)\times\cat{a}_n(a_0,a_n)
\end{equation}
for each object $b$ of $\cat{b}$; the face and degeneracy maps
are described as in \cite{hovo:mts} by composition (or evaluation) of
morphisms and by the insertion of identities respectively. The
topological realisation $B(*,\cat{a},D)$ is a $\cat{b}^{op}$-diagram
in $\cat{top}$. This definitions ensure the existence of natural
equivalences
\begin{equation}\label{bars}
B_\bullet(*,\cat{a},D)\times_{\scat{b}}E\cong
B_\bullet(*,\cat{a},D\times_{\scat{b}}E)\sands
B(*,\cat{a},D)\times_{\scat{b}}E\cong B(*,\cat{a},D\times_{\scat{b}}E)
\end{equation}
of $\cat{c}^{op}$-diagrams in $\fcat{\del$^{op}$}{top}$ and
$\cat{top}$ respectively.

\begin{exas}\label{defhoco}
If $\cat{b}=\cat{id}$, the homotopy colimit \cite{boka:hlc} of
a diagram $D\colon\cat{a}\rightarrow\cat{top}$ is given by
\[
\hoc D=B(*,\cat{a},D),
\]
as explained in \cite{hovo:mts}; using \eqref{asscot} and
\eqref{bars}, it is homeomorphic to both of
\[
B(*,\cat{a},A)\times_{\scat{a}}D\cong D
\times_{\scat{a}^{op}}B(*,\cat{a},A).
\]
In particular, $B_\bullet(*,\cat{a},*)$ is the nerve \cite{se:css}
$B_\bullet\cat{a}$ of $\cat{a}$, whose realisation is the classifying
space $B\cat{a}$ of $\cat{a}$. The natural projection $\hoc
D\rightarrow\col D$ is given by the map
\[
D\times_{\scat{a}^{op}}B(*,\cat{a},A)\longrightarrow
D\times_{\scat{a}^{op}}*,
\]
induced by collapsing $B(*,\cat{a},A)$ onto $*$.

If $\cat{a}=\cat{c}(M)$, where $M$ is an arbitrary topological monoid,
then $D$ is a left $M$-space and $B(*,\cat{c}(M),C(M))$ is a universal
contractible right $M$-space $EM$ \cite{dola:hsp}. So
\[
\hoc D=B(*,\cat{c}(M),C(M))\times_{\scat{c}(M)}D
\]
is a model for the Borel construction $EM\times_MD$.
\end{exas}

\section{Basic Constructions}\label{baco}

We choose a universal set $V$ of vertices $v_1$,\dots,$v_m$, and let $K$
denote a simplicial complex with faces $\sigma\subseteq V$. The integer
$|\sigma|-1$ is the {\it dimension\/} of $\sigma$, and the greatest such
integer is the dimension of $K$. For each $1\leq j\leq m$, the faces of
dimension less than or equal to $j$ form a subcomplex $K^{(j)}$, known
as the {\it $j$-skeleton\/} of $K$; in particular, the $1$-skeleton
$K^{(1)}$ is a graph. We abuse notation by writing $V$ for the
zero-skeleton of $K$, more properly described as $\{\{v_j\}:1\leq j\leq
m\}$. At the other extreme we have the $(m-1)$-simplex, which is the
complex containing {\it all\/} subsets of $V$; it is denoted by $2^V$
in the abstract setting and by $\varDelta(V)$ when emphasising its
geometrical realisation. Any simplicial complex $K$ therefore lies in a
chain
\begin{equation}\label{comch}
V\longrightarrow K\longrightarrow 2^V
\end{equation}
of subcomplexes. Every face $\sigma$ may also be interpreted as a
subcomplex of $K$, and so masquerades as a $(|\sigma|-1)$-simplex.

A subset $W\subseteq V$ is a {\it missing face\/} of $K$ if every
proper subset lies in $K$, yet $W$ itself does not; its dimension is
$|W|-1$. We refer to $K$ as a {\it flag complex\/}, or write that {\it
$K$ is flag}, when every missing face has two vertices. The boundary
of a planar $m$-gon is therefore flag whenever $m\geq 4$, as is the
barycentric subdivision $K'$ of an arbitrary complex $K$.  The {\it
flagification} $\Fl(K)$ of $K$ is the minimal flag complex containing
$K$ as a subcomplex, and is obtained from $K$ by adjoining every
missing face containing three or more vertices.
\begin{exa}\label{flags}
For any $n>2$, the simplest non-flag complex on $n$ vertices is the
boundary of an $(n-1)$-simplex, denoted by $\partial(n)$; then
$\Fl(\partial(n))$ is $\varDelta(n-1)$ itself.
\end{exa}
Given a subcomplex $K\subseteq L$ on vertices $V$, it is useful to
define $W\subseteq V$ as a missing face of the pair $(L,K)$ whenever
$W$ fails to lie in $K$, yet every proper subset lies in $L$.

Every finite simplicial complex $K$ gives rise to a finite category
$\cat{cat}(K)$, whose objects are the faces $\sigma$ and morphisms the
inclusions $\sigma\subseteq\tau$. The empty face $\varnothing$ is an
initial object. For any subcomplex $K\subseteq L$, the category
$\cat{cat}(K)$ is a full subcategory of $\cat{cat}(L)$; in particular,
\eqref{comch} determines a chain of subcategories
\begin{equation}\label{catch}
\cat{id}_\varnothing(m)\longrightarrow\cat{cat}(K)
\longrightarrow\cat{cat}(2^V).
\end{equation}
For each face $\sigma$, we define the {\it undercategory}
$\sigma\under\cat{cat}(K)$ by restricting attention to those objects
$\tau$ for which $\sigma\subseteq\tau$; thus $\sigma$ is an initial
object. Insisting that the inclusion $\sigma\subset\tau$ be strict
yields the subcategory $\sigma\!\Downarrow\!\cat{cat}(K)$, obtained
by deleting $\sigma$. The {\it overcategories}
$\cat{cat}(K)\under\sigma$ and $\cat{cat}(K)\!\Downarrow\!\sigma$
are defined likewise.

A complex $K$ also determines a simplicial set $S(K)$, whose
nondegenerate simplices are exactly the faces of $K$ \cite{ma:soa}.
So the nerve $B_\bullet\cat{cat}(K)$ coincides with the simplicial
set $S(\Con(K'))$, where $\Con(K')$ denotes the cone on the
barycentric subdivision of $K$, and the cone point corresponds to
$\varnothing$. More generally, $B(\sigma\under\cat{cat}(K))$ is the
cone on $B(\sigma\!\Downarrow\!\cat{cat}(K))$.

\begin{exas}\label{cubics}
If $K=V$, then $B\cat{id}_\varnothing(m)$ is the cone on $m$
disjoint points. If $K=2^V$, then $B\cat{cat}(2^V)$ is homeomorphic
to the unit cube $I^V\subset\bR^V$, and defines its canonical
simplicial subdivision; the homeomorphism maps each vertex
$\sigma\subseteq V$ to its characteristic function $\chi_\sigma$,
and extends by linearity. If $K$ is the subcomplex $\partial(m)$,
then $B\cat{cat}(\partial(m))$ is obtained from the boundary
$\partial I^m$ by deleting all faces which contain the maximal
vertex $(1,\dots,1)$.
\end{exas}

The undercategories define a $\cat{cat}(K)^{op}$-diagram
$\under\cat{cat}(K)$ in the category of small categories. It takes
the value $\sigma\under\cat{cat}(K)$ on each face $\sigma$, and the
inclusion functor
$\tau\under\cat{cat}(K)\subseteq\sigma\under\cat{cat}(K)$ on each
reverse inclusion $\tau\supseteq\sigma$. The formation of
classifying spaces yields a $\cat{cat}(K)^{op}$-diagram
$B(\vble\under\cat{cat}(K))$ in $\cat{top}_+$, which consists of
cones and their inclusions. It takes the value
$B(\sigma\under\cat{cat}(K))$ on $\sigma$ and
$B(\tau\under\cat{cat}(K))\subseteq B(\sigma\under\cat{cat}(K)$ on
$\tau\supseteq\sigma$, and its colimit is the final space
$B\cat{cat}(K)$. Following \cite{hovo:mts}, we note the isomorphism
\begin{equation}\label{catkund}
B(\vble\under\cat{cat}(K))\cong B(*,\cat{cat}(K),{\mbox{\it
CAT}}(K))
\end{equation}
of $\cat{cat}(K)^{op}$-diagrams in $\cat{top}_+$.

We refer to the cones $B(\sigma\under\cat{cat}(K))$ as {\it faces\/}
of $B\cat{cat}(K)$, amongst which we distinguish the {\it facets}
$B(v\under\cat{cat}(K))$, defined by the vertices $v$. The facets
determine the faces, according to the expression
\[
B(\sigma\under\cat{cat}(K))\;=\;
\bigcap_{v\in\sigma}B(v\under\cat{cat}(K))
\]
for each $\sigma\in K$, and form a panel structure on $B\cat{cat}(K)$
as described by Davis \cite{da:ggr}. This terminology is motivated by
our next example, which lies at the heart of recent developments in
the theory of toric manifolds.
\begin{exa}\label{polytope}
The boundary of a simplicial polytope $P$ is a simplicial complex
$K_P$, with faces $\sigma$. The {\it polar} $P^*$ of $P$ is a simple
polytope of the same dimension, whose faces $F_\sigma$ are dual to
those of $P$ (it is convenient to consider $F_\varnothing$ as $P^*$
itself). There is a homeomorphism $B\cat{cat}(K_P)\rightarrow P^*$,
which maps each vertex $\sigma$ to the barycentre of $F_\sigma$, and
transforms each face $B(\sigma\under\cat{cat}(K))$ homeomorphically
onto $F_\sigma$.
\end{exa}

Classifying the categories and functors of \eqref{catch} yields the
chain of subspaces
\begin{equation}\label{classch}
\Con(V)\longrightarrow B\cat{cat}(K)\longrightarrow I^m.
\end{equation}
So $B\cat{cat}(K)$ contains the unit axes, and is a subcomplex of
$I^m$. It is therefore endowed with the induced cubical structure,
as are all subspaces $B(\sigma\under\cat{cat}(K))$. In particular,
the simple polytope $P^*$ of Example \ref{polytope} admits a natural
cubical decomposition.

In our algebraic context, we utilise the category $\cat{grp}$ of
discrete groups and homomorphisms. Many constructions in $\cat{grp}$ may
be obtained by restriction from those we describe in $\cat{tmon}$, and
we leave readers to provide the details. In particular, $\cat{grp}$ is a
full subcategory of $\cat{tmg}$, and is $\cat{top}_+$-complete and
-cocomplete.

Given a commutative ring $Q$ (usually the integers, or their
reduction mod 2), we consider the category ${}_Q\cat{mod}$ of left
$Q$-modules and $Q$-linear maps, which is symmetric monoidal with
respect to the tensor product $\otimes_Q$ and closed under
$(Z,Y)\mapsto{}_Q\cat{mod}(Y,Z)$. We usually work in the related
category $\cat{g$_Q$mod}$ of connected graded modules of finite
type, or more particularly in the categories $\cat{g$_Q$calg}$ and
$\cat{g$_Q$cocoa}$, which are dual; the former consists of augmented
commutative $Q$-algebras and their homomorphisms, and the latter of
supplemented cocommutative $Q$-coalgebras and their coalgebra maps.

As an object of ${}_Q\cat{mod}$, the polynomial algebra $Q[V]$ on
$V$ has a basis of monomials $v_W=\prod_Wv_j$, for each multiset $W$
on $V$. Henceforth, we assign a common dimension $d(v_j)>0$ to the
vertices $v_j$ for all $1\leq j\leq m$, and interpret $Q[V]$ as an
object of $\cat{g$_Q$calg}$; products are invested with appropriate
signs if $d(v_j)$ is odd and $2Q\neq0$. Then the quotient map
\[
Q[V]\longrightarrow Q[V]/(v_\lambda:\lambda\notin K)
\]
is a morphism in $\cat{g$_Q$calg}$, whose target is known as the graded
{\it Stanley-Reisner $Q$-algebra} of the simplicial complex $K$, and
written $\SR_Q(K)$. This ring is a fascinating invariant of $K$, and
reflects many of its combinatorial and geometrical properties, as
explained in \cite{st:cca}. Its $Q$-dual is a graded incidence coalgebra
\cite{joro:cbc}, which we denote by $\SR^Q(K)$.

We define a $\cat{cat}(K)^{op}$-diagram $D_K$ in $\cat{top}_+$ as
follows. The value of $D_K$ on each face $\sigma$ is the discrete space
$\sigma_+$, obtained by adjoining $+$ to the vertices, and the value on
$\tau\supseteq\sigma$ is the projection $\tau_+\rightarrow\sigma_+$,
which fixes the vertices of $\sigma$ and maps the vertices of
$\tau\setminus\sigma$ to $+$.
\begin{defn}\label{expodefs}
Given objects $(X,*)$ of $\cat{top}_+$ and $M$ of $\cat{tmg}$, the
{\it exponential diagrams} $X^K$ and $M^K$ are the cotensor homsets
$\Hom_{\cat{id}}(D_K,X)$ and $\Hom_{\cat{id}}(D_K,M)$ respectively;
they are $\cat{cat}(K)$-diagrams in $\cat{top}_+$ and $\cat{tmg}$.
Alternatively, they are the respective compositions of the
exponentiation functors
$X^{(\vble)}\colon\cat{top}_+^{op}\rightarrow\cat{top}_+$ and
$M^{(\vble)}\colon\cat{top}_+^{op}\rightarrow\cat{tmg}$ with
$D_K^{op}$.
\qed
\end{defn}
So the value of $X^K$ on each face $\sigma$ is the product space
$X^\sigma$, whose elements are functions $f\colon\sigma\rightarrow X$,
and the value of $X^K$ on $\sigma\subseteq\tau$ is the inclusion
$X^\sigma\subseteq X^\tau$ obtained by extending $f$ over $\tau$ by the
constant map $*$. The space $X^\varnothing$ consists only of $*$. In the
case of $M^K$, each $M^\sigma$ is invested with pointwise
multiplication, so $H^K$ takes values in $\cat{grp}$ for a discrete
group $H$.

In $\cat{g$_Q$calg}$, we define a $\cat{cat}(K)^{op}$-diagram $Q[K]$
by analogy. Its value on $\sigma$ is the graded polynomial algebra
$Q[\sigma]$, and on $\tau\supseteq\sigma$ is the projection
$Q[\tau]\rightarrow Q[\sigma]$. We denote the dual
$\cat{cat}(K)$-diagram $\Hom_{\cat{id}}(Q[K],Q)$ by $Q\langle
K\rangle$, and note that it lies in $\cat{g$_Q$cocoa}$. Its value on
$\sigma$ is the free $Q$-module $Q\langle S(\sigma)\rangle$ generated
by simplices $z$ in $S(\sigma)$, and on $\sigma\subseteq\tau$ is the
corresponding inclusion of coalgebras. The coproduct is given by
$\delta(z)=\sum z_1\otimes z_2$, where the sum ranges over all
partitions of $z$ into subsimplices $z_1$ and $z_2$.

When $Q=\bZ/2$ we let the vertices have dimension $1$. Every
monomial $v_U$ therefore has dimension $|U|$ in the graded algebra
$\bZ/2[\sigma]$, and every $j$-simplex in $S(\sigma)$ has dimension
$j+1$ in $\bZ/2\langle S(\sigma)\rangle$. We refer to this as the
{\it real case}. When $Q=\bZ$ we consider two possibilities. First
is the {\it complex case}, in which the vertices have dimension $2$,
so that the additive generators of $\bZ[\sigma]$ and $\bZ\langle
S(\sigma)\rangle$ have twice the dimension of their real
counterparts. Second is the {\it exterior case}, in which the
dimension of the vertices reverts to $1$. Every squarefree monomial
$v_U$ then has dimension $|U|$ in $\bZ[\sigma]$, and
anticommutativity ensures that every monomial containing a square is
zero; every $j$-face of $\sigma$ has dimension $j+1$ in $\bZ\langle
S(\sigma)\rangle$, and every degenerate $j$-simplex $z$ represents
zero. To distinguish between the complex and exterior cases, we
write $Q$ as $\bZ$ and $\wedge$ respectively.

In the real and complex cases, \daaja\ \cite{daja:cpc} introduce
homotopy types $\djs_\bR(K)$ and $\djs_\bC(K)$. The cohomology rings
$H^*(\djs_\bR(K);\bZ/2)$ and $H^*(\djs_\bC(K);\bZ)$ are isomorphic to
the graded Stanley-Reisner algebras $\SR_{\bZ/2}(K)$ and $\SR_\bZ(K)$
respectively. We shall deal with the exterior case below, and discuss
alternative constructions for all three cases. We write $\djs(K)$ as a
generic symbol for \daaja's homotopy types, and refer to them as {\it
\daja\/} spaces for $K$. They are represented by objects in \cat{top}.

\section{Colimits}\label{co}

In this section we introduce the colimits which form our main topic of
discussion, appealing to the completeness and cocompleteness of
$\cat{t}$ and $\cat{tmg}$ as described in Section \ref{capr}. We
consider colimits of the diagrams $X^K$, $M^K$, $G^K$, and $Q\langle
K\rangle$ in the appropriate categories, and label them $\col^+X^K$,
$\coltmg M^K$, $\coltmg G^K$, and $\col Q\langle K\rangle$
respectively. Similarly, we write the limit of $Q[K]$ as $\lim Q[K]$. As
we shall see, these limits and colimits coincide with familiar
constructions in several special cases.

As an exercise in acclimatisation, we begin with the diagrams
associated to \eqref{catch}. Exponentiating with respect to $(X,*)$
and taking colimits provides the chain of subspaces
\begin{equation}\label{ptopch}
\bigvee_{j=1}^mX_j\longrightarrow \col^+X^K\longrightarrow X^m,
\end{equation}
thereby sandwiching $\col^+X^K$ between the axes and the cartesian
power. On the other hand, using an object $M$ of \cat{tmg} yields
the chain of epimorphisms
\begin{equation}\label{tmonch}
\bigfree_{j=1}^mM_j\longrightarrow\coltmg M^K\longrightarrow M^m,
\end{equation}
giving a presentation of $\coltmg M^K$ which lies between the
$m$-fold free product of $M$ and the cartesian power.

The following example emphasises the influence of the underlying
category on the formation of colimits, and is important later.
\begin{exa}\label{fatw}
If $K$ is the non-flag complex $\partial(m)$ of Example {\rm
\ref{flags}} (where $m>2$), then $\col^+X^K$ is the {\rm fat wedge\/}
subspace $\{(x_1,\dots,x_m):x_j=*\;\text{for some $1\leq j\leq m$}\}$;
on the other hand, $\coltmg M^K$ is isomorphic to $M^m$ itself.
\end{exa}

By construction, $\coltmg C_2^K$ in \cat{grp} enjoys the presentation
\[
\langle a_1,\dots,a_m : a_j^2=1, (a_ia_j)^2=1\;\text{for all
$\{v_i,v_j\}$ in $K$}\rangle
\]
and is isomorphic to the {\it right-angled Coxeter group\/}
$\Cox(K^{(1)})$ determined by the 1-skeleton of $K$. Readers should
not confuse $K^{(1)}$ with the more familiar Coxeter graph of the
group, which is almost its complement!

Similarly, $\coltmg C^K$ has the presentation
\[
\langle b_1,\dots,b_m :
[b_i,b_j]=1\;\text{for all $\{v_i,v_j\}$ in $K$}\rangle
\]
(where $[b_i,b_j]$ denotes the commutator $b_ib_jb_i^{-1}b_j^{-1}$),
and so is isomorphic to the {\it right-angled Artin group\/}
$\Art(K^{(1)})$. Such groups are sometimes called {\it graph groups},
and are special examples of {\it graph products} \cite{ch:ecg}. As
explained to us by Dave Benson, neither should be confused with the
{\it graphs of groups\/} described in \cite{wa:tmg}.

In the continuous case, we define the {\it circulation group\/}
$\Cir(K^{(1)})$ as $\coltmg T^K$ in $\cat{tmg}$. Every element of $\Cir
(K^{(1)})$ may therefore be represented as a word
\begin{equation}\label{cirword}
t_{i_1}(1)\cdots t_{i_k}(k),
\end{equation}
where $t_{i_j}(j)$ lies in the $i_j$th factor $T_{i_j}$ for each $1\leq
j\leq k$. Two elements $t_r\in T_r$ and $t_s \in T_s$ commute whenever
$\{r,s\}$ is an edge of $K$.

We shall use $G$ as a generic symbol for any one of the groups $C_2$,
$C$, or $T$.

Following \eqref{tmonch}, we abbreviate the generating subgroups
$G^{v_j}<\col G^K$ to $G_j$, where $1\leq j\leq m$, and call them the
{\it vertex groups}. Since $\coltmg G^K$ is presented as a quotient of
the free product $\bigfree_{j=1}^mG_j$, its elements $g$ may be
assigned a wordlength $l(g)$. In addition, the arguments of
\cite{cafo:pcc} apply to decompose every $g$ {\it from the right\/} as
\begin{equation}\label{wdecomp}
g\;=\;\prod_{j=1}^ns_j(g)
\end{equation}
for some $n\leq l(g)$, where each subword $s_j(g)$ contains the
maximum possible number of mutually commuting letters, and is unique.

Given any subset $W\subseteq V$ of vertices, we write $K_W$ for the
complex obtained by restricting $K$ to $W$. The following Lemma is a
simple restatement of the basic properties of $\coltmg G^K$.
\begin{lem}\label{awsubs}
We have that
\begin{enumerate}
\item
the subgroup $\coltmg G^{K_W}\leq\coltmg G^K$ is abelian if and only if
$K_W^{(1)}$ is a complete graph, in which case it is isomorphic to
$G^W$;
\item
when $K$ is flag, each subword $s_j(g)$ of \eqref{wdecomp} lies in a
subgroup $G^{\sigma_j}$ for some face $\sigma_j$ of $K$.
\end{enumerate}
\end{lem}

Other algebraic examples of our colimits relate to the Stanley-Reisner
algebras and coalgebras of $K$. By construction, there are algebra
isomorphisms
\begin{equation}\label{sralgisos}
\lim\bZ/2[K]\cong\SR_{\bZ/2}(K),\quad
\lim\bZ[K]\cong\SR_\bZ(K),\sands
\lim\wedge[K]\cong\SR_\wedge(K),
\end{equation}
where the limits are taken in $\cat{g$_\bZ$calg}$. Dually, there are
coalgebra isomorphisms
\begin{equation}\label{srcoisos}
\col\bZ/2\langle K\rangle\cong\SR^{\bZ/2}(K),\quad
\col\bZ\langle K\rangle\cong\SR^\bZ(K),\sands
\col\wedge\langle K\rangle\cong\SR^\wedge(K)
\end{equation}
in $\cat{g$_\bZ$cocoa}$. The analogues of \ref{ptopch} display these
limits and colimits as
\begin{equation}\label{analog}
\bigoplus_{j=1}^m Q[v_j]\longleftarrow\lim Q[K]\longleftarrow Q[V]
\sands
\bigoplus_{j=1}^m\DP^Q(v_j)\longrightarrow\col
Q\langle K\rangle\longrightarrow\DP^Q(V)
\end{equation}
respectively; here $\DP^Q(W)$ denotes the divided power $Q$-coalgebra
of multisets on $W\subseteq V$, graded by dimension.

If we let $(X,*)$ be one of the pairs $(BC_2,*)$, $(BT,*)$, or $(BC,*)$,
then simple arguments with cellular chain complexes show that the
cohomology rings $H^*(\col^+(BC_2)^K;\bZ/2)$, $H^*(\col^+(BT)^K;\bZ)$,
and $H^*(\col^+(BC)^K;\bZ)$ are isomorphic to the limits
\eqref{sralgisos} respectively. Similarly, the homology coalgebras are
isomorphic to the dual coalgebras \eqref{srcoisos}. In cohomology, these
observations are due to Buchstaber and Panov \cite{bupa:tact} in the
real and complex cases, and to Kim and Roush \cite{kiro:hca} in the
exterior case (at least when $K$ is $1$-dimensional). In homology, they
may be made in the context of incidence coalgebras, following
\cite{rasc:cmc}. In both cases, the maps of \eqref{ptopch} induce the
homomorphisms \eqref{analog}.

Such calculations do not themselves identify $\col^+(BC_2)^K$ and
$\col^+(BT)^K$ with \daaja's constructions.  Nevertheless, Buchstaber
and Panov provide homotopy equivalences
$\col^+(BC_2)^K\simeq\djs_\bR(K)$ and $\col^+(BT)^K\simeq\djs_\bC(K)$,
which also follow from Corollary \ref{hofibkl} below; the Lemma
yields a corresponding equivalence in the exterior case. Of course,
$\col^+(BC)^K$ is a subcomplex of the $m$-dimensional torus $(S^1)^m$,
and is therefore finite.

In due course, we shall use these remarks to interpret the following
proposition in terms of \daja\ spaces. The proof for $G=C_2$ is
implicit in \cite{daja:cpc}, and for $G=C$ is due to Kim and Roush
\cite{kiro:hca}.
\begin{prop}\label{djkrthm}
When $G=C_2$ or $C$, there is a homotopy equivalence
\[
\col^+(BG)^K\;\simeq\;B\coltmg G^K
\]
for any flag complex $K$.
\end{prop}

Since both cases are discrete, $B\coltmg G^K$ is, of course, an
Eilenberg-Mac Lane space; Charney and Davis \cite{chda:fka} have since
identified good models for $BA$, given {\it any\/} Artin group
$A$. Proposition \ref{djkrthm} fails for arbitrary complexes $K$, as
our next examples show.
\begin{exas}\label{exaso}
Proposition {\rm \ref{djkrthm}} applies when $K=V$, because the discrete
complex is flag; then $\coltmg G^K$ is isomorphic to the free product of
$m$ copies of $G$, whose classifying space is the $m$-fold wedge
$\bigvee_{j=1}^m BG_j$ (by \cite{brha:tgo}, for example). On the other
hand, when $K$ is the non-flag complex $\partial(m)$, Example {\rm
\ref{fatw}} confirms that $B\coltmg G^K$ is $BG^m$, whereas
$\col^+(BG)^K$ is the fat wedge subspace.
\end{exas}

These examples apply unchanged to the case $G=T$, and serve to motivate
our extension of Proposition \ref{djkrthm} to the complex case in
Proposition \ref{hkkcontract} below. So far as $C_2$ and $C$ are
concerned, the Proposition asserts that certain homotopy homomorphisms
\begin{equation}\label{hoho}
h_K\colon\varOmega\col^+(BG)^K\longrightarrow\coltmg G^K
\end{equation}
are homotopy equivalences when $K$ is flag. We therefore view the
$h_K$ as modelling the loop spaces; in the complex case, they express
$\varOmega\col^+(BT)^K$ in terms of the circulation groups $\coltmg
T^K$. In Section \ref{hocotomo} we will use homotopy colimits to
describe analogues of $h_K$ for all complexes $K$.

Our interest in the loop spaces $\varOmega\col^+(BG)^K$ has been
stimulated by several ongoing programmes in combinatorial algebra. For
example, Herzog, Reiner, and Welker \cite{herewe:cli} discuss
combinatorial issues associated with calculating the $k$-vector spaces
$\Tor^{S\hspace{-.3mm}R_k(K)}(k,k)$ over an arbitrary ground field $k$,
and refer to \cite{gule:hlr} for historical background. Such
calculations have applications to diagonal subspace arrangements, as
explained by Peeva, Reiner and Welker \cite{perewe:crd}. Since these
$\Tor$ spaces also represent the $E_2$-term of the Eilenberg-Moore
spectral sequence for $H^*(\varOmega\djs(K);k)$, it seems well worth
pursuing geometrical connections. We consider the algebraic implications
elsewhere \cite{para:hhc}.

\section{Fibrations and homotopy colimits}\label{fihoco}

In this section we apply the theory of homotopy colimits to study
various relevant fibrations and their geometrical interpretations.
Some of the results appear in \cite{bupa:tact}, but we believe that
our approach offers an attractive and efficient alternative, and eases
generalisation. We refer to \cite{hovo:mts} and \cite{vo:hlc} for the
notation and fundamental properties of homotopy colimits. Several of
the results we use are also summarised in \cite{wezizi:hcc}, together
with additional information on combinatorial applications.

We begin with a general construction, based on a well-pointed
topological group $\varGamma$ and a diagram
$H\colon\cat{a}\rightarrow\cat{tmg}$ of closed subgroups and their
inclusions. We assume that the maps of the classifying diagram
$BH\colon\cat{a}\rightarrow\cat{top}_+$ are cofibrations, and that
the Projection Lemma \cite{wezizi:hcc} applies to the natural
projection $\hoc^+ BH\rightarrow\col^+BH$, which is therefore a
homotopy equivalence. The cofibrations $BH(a)\rightarrow B\varGamma$
correspond to the canonical map $f_H\colon\col^+BH\rightarrow
B\varGamma$ under the homeomorphism \eqref{eneqcoeq}.

By Examples \ref{simco} the coset spaces $\varGamma/H(a)$ define an
$\cat{a}\times\cat{c}(\varGamma)$ diagram $\varGamma/H$ in $\cat{top}$,
and by Examples \ref{defhoco} the cofibration $BH(a)\rightarrow
B\varGamma$ is equivalent to the fibration
\[
B\big(*,\cat{c}(\varGamma),C(\varGamma)
\times_{\scat{c}(\varGamma)}\varGamma/H(a)\big)\longrightarrow
B\cat{c}(\varGamma)
\]
for each object $a$ of $\cat{a}$. So $f_H$ is equivalent to
\[
\hoc^+B(*,\cat{c}(\varGamma),C(\varGamma)
\times_{\scat{c}(\varGamma)}\varGamma/H)\longrightarrow B\varGamma
\]
in the homotopy category of spaces over $B\varGamma$, where the homotopy
colimit is taken over $\cat{a}$.
\begin{prop}\label{fibcofib}
The homotopy fibre of $f_H$ is the homotopy colimit $\hoc^+\varGamma/H$.
\end{prop}
\begin{proof}
We wish to identify the homotopy fibre of the projection
\[
B\big(*,\cat{a},B(*,\cat{c}(\varGamma),C(\varGamma)
\times_{\scat{c}(\varGamma)}\varGamma/H)\big)\longrightarrow B\varGamma.
\]
But we may rewrite the total space as
$B(*,\cat{a},\varGamma/H)\times_{\scat{c}(\varGamma)^{op}}
B(*,\cat{c}(\varGamma),C(\varGamma))$, and therefore as
$B(*,\cat{c}(\varGamma),C(\varGamma))
\times_{\scat{c}(\varGamma)}B(*,\cat{a},\varGamma/H)$, using
\eqref{bars} and Examples \ref{defhoco}. So the homotopy fibre is
$B(*,\cat{a},\varGamma/H)$, as required.
\end{proof}

Given a pair of simplicial complexes $(L,K)$ on vertices $V$, we let
$\cat{a}=\cat{cat}(K)$, and choose $\varGamma=\coltmg G^L$ and $H=G^K$;
we also abbreviate the diagram $\varGamma/H$ to $L/K$. Then $f_H$ is the
induced map
\begin{equation}
f_{K,L}\colon\col^+(BG)^K\longrightarrow B\coltmg G^L,
\end{equation}
and the Projection Lemma applies to $(BG)^K$ because the maps
$\col^+(BG)^{K\Downarrow\sigma}\rightarrow BG^\sigma$ are closed
cofibrations for each face $\sigma$. So we have the following
corollary to Proposition \ref{fibcofib}.
\begin{cor}\label{hofibkl}
The homotopy fibre of $f_{K,L}$ is the homotopy colimit $\hoc L/K$,
and is homeomorphic to the identification space
\begin{equation}\label{identsp}
\left(B\cat{cat}(K)\times\coltmg G^L\right)/\sim,
\end{equation}
where $(p,gh)\sim(p,g)$ whenever $h\in G^\sigma$ and $p$ lies in the
face $B(\sigma\under\cat{cat}(K))$.
\end{cor}
\begin{proof}
By \eqref{catkund}, the homotopy colimit $B(*,\cat{cat}(K),L/K)$ may
be expressed as
\[
B(\vble\under\cat{cat}(K))\times_{\scat{cat}(K)}L/K,
\]
and the inclusions $B(\sigma\under\cat{cat}(K))\subseteq B\cat{cat}(K)$
induce a homeomorphism with \eqref{identsp}.
\end{proof}
For future use, we write $\mu$ for the canonical action of
$\coltmg G^L$ on $B(*,\cat{cat}(K),L/K)$.

We note that $f_{K,L}$ coincides with the right-hand map of
\eqref{ptopch} when $L=2^V$ and $X=BG$; the cases in which $K=L$
(abbreviated to $f_K$) and $L=\Fl(K)$ also feature below. The space
$\hoc 2^V/K$ plays a significant r\^ole in \cite{daja:cpc}, where it
is described as the identification space of Corollary \ref{hofibkl}
and denoted by $\mathcal{Z}_P$ (with $P$ the dual of $K$, in the sense
of Example \ref{polytope}). To emphasise this connection, we write
$\hoc L/K$ as $\mathcal{Z}_G(K,L)$, which we abbreviate to
$\mathcal{Z}_G(K)$ when $K=L$. It appears repeatedly below, by virtue
of Proposition \ref{fibcofib}. Our examples assume that $L=2^V$, and
continue the theme of Examples \ref{exaso}.
\begin{exas}\label{exast}
If $K=V$ then $\mathcal{Z}_G(K,2^V)$ is the homotopy fibre of
$\bigvee_{j=1}^m BG_j\rightarrow BG^m$, the inclusion of the axes; it
has been of interest to homotopy theorists for many years. If $K$ is
the non-flag complex $\partial(m)$, then $\mathcal{Z}_G(K,2^V)$ is
homotopy equivalent to $S^{m-1}$ for $G=\bZ/2$, and $S^{2m-1}$ for
$G=T$.
\end{exas}
The second of these examples may be understood by noting that the
inclusion of the fat wedge in $BG^m$ has the Thom complex of the
external product $\zeta^m$ of Hopf bundles as its cofibre.

\daaja\ \cite{daja:cpc} prove that the mod 2 cohomology ring of
$EC_2^m\times_{C_2^m}\mathcal{Z}_{C_2}(K,2^V)$ and the integral
cohomology ring of $ET^m\times_{T^m}\mathcal{Z}_T(K,2^V)$ are
isomorphic to the Stanley Reisner algebras $\SR_{\bZ/2}^*(K)$ and
$\SR_\bZ^*(K)$ respectively. In view of Corollary \ref{hofibkl} (in
the case $L=2^V$), we regard the spaces $\col^+(BG)^K$ and the \daja\
homotopy types as interchangeable from this point on.

The canonical projection $\mathcal{Z}_G(K,L)\rightarrow
B\cat{cat}(K)$ is obtained by factoring out the action
$\mu$ of $\coltmg G^L$ on $\hoc L/K$. The cubical structure
\eqref{classch} of the quotient lifts to an associated decomposition
of $\mathcal{Z}_G(K,L)$; when $G=T$ and $L=2^V$, for example, we
recover the description of \cite{bupa:tact} and \cite{daja:cpc} in
terms of polydiscs and tori.

The action $\mu$ has other important properties.
\begin{prop}\label{isots}
The isotropy subgroups of $\mu$ are the conjugates $wG^\sigma
w^{-1}<\coltmg G^L$, where $\sigma$ ranges over the faces of $K$.
\end{prop}
\begin{proof}
It suffices to note from Corollary \ref{hofibkl} that each point
$[x,wG^\sigma]$ is fixed by $wG^\sigma w^{-1}<\coltmg G^L$, for any
$x\in B(\sigma\under\cat{cat}(K))$.
\end{proof}
\begin{cor}\label{commsub}
The commutator subgroup of $\coltmg G^L$ acts freely on
$\mathcal{Z}_G(K,L)$ under $\mu$.
\end{cor}
\begin{proof}
The isotropy subgroups are abelian, and so have trivial intersection
with the commutator subgroup.
\end{proof}

When $K=L$ and $G=C_2$, Proposition \ref{isots} strikes a familiar
chord. The {\it parabolic\/} subgroups of a Coxeter group $H$ are the
conjugates $w\varGamma w^{-1}$ of certain subgroups $\varGamma$,
generated by subsets of the defining Coxeter system; when $H$ is
right-angled, and therefore takes the form $\Cox(K^{(1)})$, such
subgroups are abelian. When $L=2^V$, each subgroup $wG^\sigma w^{-1}$
reduces to $G^\sigma$.  In this case, Proposition \ref{isots} implies
that the isotropy subgroups form an exponential
$\cat{cat}^{op}(K)$-diagram in \cat{tgrp}, which assigns
$G^\sigma$ to the face $\sigma$ and the quotient homomorphism
$G^\tau\rightarrow G^\sigma$ to the reverse inclusion
$\tau\supseteq\sigma$.

As detailed in \cite{bupa:tact}, the homotopy fibre
$\mathcal{Z}_G(K,2^V)$ is closely related to the theory of subspace
arrangements and their auxiliary spaces. These spaces are defined in
each of the real, complex, and exterior cases, and will feature below;
we introduce them here as homotopy colimits.

Given a pointed space $(Y,0)$, we let $Y_\times$ denote $Y\setminus
0$. For any subset $W\subseteq V$, we write $Y_W\subseteq Y^V$ for the
{\it coordinate subspace\/} of functions $f\colon V\rightarrow Y$ for
which $f(W)=0$. The set of subspaces
\[
\mathcal{A}_Y(K) = \{Y_W:W\notin K\}
\]
is the associated {\it arrangement\/} of $K$, whose {\it complement\/}
$U_Y(K)$ is given by the equivalent formulae
\begin{equation}\label{altcomp}
Y^V\setminus{\textstyle\bigcup_{W\notin K}}Y_W
\;=
\;\{f:f^{-1}(0)\in K\}.
\end{equation}
The $\cat{cat}(K)$-diagram $Y(K)$ associates the function
space $Y(\sigma)=\{f:f^{-1}(0)\subseteq\sigma\}$ to each face $\sigma$,
and the inclusion $Y(\sigma)\subseteq Y(\tau)$ to each morphism
$\sigma\subseteq\tau$. It follows that $Y(\sigma)$ is homeomorphic to
$Y^\sigma\times(Y_\times^{V\setminus\sigma})$, and that $U_Y(K)$ is
$\col Y(K)$.

The exponential $\cat{cat}(K)$-diagram
$Y_\times^{V\setminus K}$ associates $Y_\times^{V\setminus\sigma}$ to
$\sigma$; when $Y$ is contractible, we may therefore follow Proposition
\ref{fibcofib} by combining the Projection Lemma and Homotopy Lemma of
\cite{wezizi:hcc} to obtain a homotopy equivalence
\begin{equation}\label{hocoarr}
\hoc Y_\times^{V\setminus K}
\;\simeq\; U_Y(K).
\end{equation}

Now let us write $\bF$ for one of the fields $\bR$ or $\bC$. The study
of the {\it coordinate subspace arrangements} $\mathcal{A}_\bF(K)$,
together with their complements, is a special case of a well-developed
theory whose history is rich and colourful (see \cite{bj:sa}, for
example). In the exterior case, we replace $\bF$ by the union of a
countably infinite collection of 1-dimensional cones in $\bR^2$, which
we call a {\it $1$-star} and write as $\bE$. So $\bE^V$ is an $m$-star;
it is homeomorphic to the union of countably many $m$-dimensional cones
in $(\bR^2)^V$, obtained by taking products.

As $G$ ranges over $C_2$, $T$ and $C$, we let $\bF$ denote $\bR$, $\bC$
and $\bE$ respectively. In all three cases, the natural inclusion of $G$
into $\bF_\times$ is a cofibration, and $\bF_\times$ retracts onto its
image. So \eqref{hocoarr} applies, and may be replaced by the
corresponding equivalence
\begin{equation}\label{hocoars}
\hoc G^{V\setminus K}\;\simeq\; U_\bF(K).
\end{equation}
\begin{prop}\label{hofili}
The space $\mathcal{Z}_G(K,2^V)$ is homotopy equivalent to $U_\bF(K)$,
for any complex $K$.
\end{prop}
\begin{proof}
Substitute $L=2^V$ in Corollary \ref{hofibkl} and apply \eqref{hocoars}.
\end{proof}

By specialising certain results of \cite{wezizi:hcc} and
\cite{zizi:hts}, we may also describe $\big(\bigcup_{W\notin
K}\bF_W\big)\setminus 0$ as a homotopy colimit. This space is dual to
$U_\bF(K)$, and appears to have a more manageable homotopy type in
many relevant cases. For $G=C_2$ and $T$, a version of Proposition
\ref{hofili} features prominently in \cite{bupa:tact}.

The following examples illustrate Proposition \ref{hofili}, in the
light of Examples \ref{exast}.
\begin{exas}\label{arrncomps}
For $m>2$ and $G=T$, the subspace arrangements of the discrete
complex $V$ and the non-flag complex $\partial(m)$ are given by
\[
\big\{\{z:z_j=z_k=0\}:1\leq j<k\leq m\big\}\sands\{0\}
\]
respectively; the corresponding complements are
\[
\{z:z_j=0\Rightarrow z_k\neq 0\}\sands \bC^m\setminus 0.
\]
The former is homotopy equivalent to a wedge of spheres, and the
latter to $S^{2m-1}$.
\end{exas}

\section{Flag complexes and connectivity}\label{flcoco}

In this section, we examine the homotopy fibre $\mathcal{Z}_G(K,L)$
more closely. The results form the basis of our model for
$\varOmega\djs(K)$ when $K$ is flag, and enable us to measure the
extent of its failure for general $K$.

We consider a flag complex $K$, and substitute $K=L$ into Corollary
\ref{hofibkl} to deduce that $\mathcal{Z}_G(K)$ is the homotopy fibre
of the cofibration $f_K\colon\djs(K)\rightarrow B\coltmg G^K$. It is
helpful to abbreviate $B(\sigma\under\cat{cat}(K))$ to $B(\sigma)$
thoughout the following argument.
\begin{prop}\label{hkkcontract}
The cofibration $f_K$ is a homotopy equivalence whenever $K$ is flag.
\end{prop}
\begin{proof}
We prove that $\mathcal{Z}_G(K)$ is contractible.

For any face $\sigma\in K$, the space $(\coltmg G^K)/G^\sigma$ inherits
an increasing filtration by subspaces $(\coltmg G^K)_i/G^\sigma$,
consisting of those cosets $wG^\sigma$ for which a representing
element satisfies $l(w)\leq i$. We may therefore define a
$\cat{cat}(K)$-diagram $K_i/K$, which assigns $(\coltmg G^K)_i/G^\sigma$
to each face $\sigma$ and the corresponding inclusion to each
inclusion $\sigma\subseteq\tau$. By construction, $\mathcal{Z}_G(K)$
is filtered by the subspaces $\hoc K_i/K$ and each inclusion $\hoc
K_{i-1}/K\subset\hoc K_i/K$ is a cofibration. We proceed by induction
on $i$.

For the base case $i=0$, we observe that $(\coltmg G^K)_0/G^\sigma$ is
the single point $eG^\sigma$ for all values of $\sigma$. Thus $\hoc
K_0/K$ is homeomorphic to $B(\varnothing)$, and is indeed contractible.
To make the inductive step, we assume that $\hoc K_i/K$ is contractible
for all $i<n$, and write $Q_n$ for the quotient space $(\hoc
K_n/K)/(\hoc K_{n-1}/K)$. It then suffices to prove that $Q_n$ is
contractible.

Every point of $Q_n$ has the form $(x,wG^\sigma)$, for some $x\in
B(\sigma)$ and some $w$ of length $n$. If the final letter of $w$ lies
in $G^\sigma$, then $(x,wG^\sigma)$ is the basepoint of
$Q_n$. Otherwise, we rewrite $w$ as $w's$ by \eqref{wdecomp}, where
$s$ contains the maximum possible number of mutually commuting
letters. These determine a subset $\chi\subseteq V$, and Lemma
\ref{awsubs} confirms that $K^{(1)}$ contains the complete graph on
vertices $\chi$. Since $K$ is flag, we deduce that $2^\chi\in K$, and
therefore that $(x,w'G^\chi)$ is the basepoint of $Q_n$. To describe a
contraction of $Q_n$, we may find a canonical path $p$ in
$\cat{cat}_\varnothing(K)$, starting at $x$ and finishing at some $x'$
in $B(\chi)$; of course $p$ must vary continuously with
$(x,wG^\sigma)$, and lift to a corresponding path in $Q_n$. If $x$ is
a vertex of $B(\sigma)$, we choose $p$ to run at constant speed along
the edge from $x$ to the cone point $\varnothing$, and again from
$\varnothing$ to the vertex $\chi\in B(\chi)$. If $x$ is an interior
point of $B(\sigma)$, we extend the construction by linearity. Then
$p$ lifts to the path through $(p(t),w)$ for all $0<t<1$, as required.
\end{proof}

Proposition \ref{hkkcontract} leads to the study of
$f_{K,L}\colon\djs(K)\rightarrow B\coltmg G^L$ for any subcomplex
$K\subseteq L$. We consider the missing faces of $K$ with three or
more vertices and write $c(K)\geq 2$ for their minimal dimension. We
let $d(K)$ denote $c(K)-1$ when $G=C_2$ or $C$, and $2c(K)$ when
$G=T$; thus $K$ is flag if and only if $c(K)$ (and therefore $d(K)$)
is infinite. Finally, we define
\[
c(K,L)=
\begin{cases}
c(K)&\text{if $L\subseteq\Fl(K)$}\\
1&\text{otherwise},
\end{cases}
\]
and let $d(K,L)$ be given by $c(K,L)-1$ or $2c(K,L)$ as before.
\begin{thm}\label{lkequiv}
For any subcomplex $K\subseteq L$, the cofibration $f_{K,L}$ is a
$d(K,L)$-equivalence.
\end{thm}
\begin{proof}
We may factorise $f_{K,L}$ as
\[
\djs(K)\longrightarrow\djs(Fl(K))\longrightarrow\djs(\Fl(L))
\longrightarrow B\coltmg G^{Fl(L)}.
\]
The first map is induced by flagification, and is a $d(K)$-equivalence
by construction. The second is the identity if $L\subseteq\Fl(K)$;
otherwise, it is $0$-connected when $G=C_2$ or $C$, and $2$-connected
when $G=T$. The third map is $f_{Fl(L)}$, and an equivalence by
Proposition \ref{hkkcontract}.
\end{proof}

Theorem \ref{lkequiv} suggests our first model for $\varOmega\djs(K)$.
\begin{prop}\label{dequiv}
There is a homotopy homomorphism
$h_K\colon\varOmega\djs(K)\rightarrow\coltmg G^K$, which is a
$(d(K)-1)$-equivalence for any complex $K$; in particular, it is an
equivalence if $K$ is flag.
\end{prop}
\begin{proof}
Applying Theorem \ref{lkequiv} with $K=L$ implies that $\varOmega
f_K\colon\varOmega\djs(K)\rightarrow\varOmega B\coltmg G^K$ is a
$(d(K)-1)$-equivalence. The result follows by composing with the
canonical homotopy homomorphism $\varOmega BH\rightarrow H$, which
exists for any topological group $H$.
\end{proof}
When $L=2^V$, the missing faces of $(2^V,K)$ are precisely the
non-faces of $K$. In this case only, we write their minimal dimension
as $c'(K)$.

It is instructive to consider the homotopy commutative diagram
\begin{equation}\label{diago}
\begin{CD}
\mathcal{Z}_G(K,L)@>id>>\mathcal{Z}_G(K,L)@>>>*\\
@VVpV@VVV@VVV\\
\mathcal{Z}_G(K,2^V)@>>>\djs(K)@>\;f_{K,2^V}\;>>BG^m\\
@VV\gamma V@VVf_{K,L}V@VVidV\\
B[G,L]@>>>B\coltmg G^L@>Ba>>BG^m
\end{CD}
\end{equation}
of fibrations, where $a$ is the abelianisation homomorphism and
$[G,L]$ denotes the commutator subgroup of $\coltmg G^L$. By Theorem
\ref{lkequiv}, $\mathcal{Z}_G(K,L)$ and $\mathcal{Z}_G(K,2^V)$ are
$(d(K,L)-1)$- and $(d'(K)-1)$-connected respectively, where
$d(K,L)\geq d'(K)$ by definition. In fact $\mathcal{Z}_G(K,2^V)$ is
$d'(K)$-connected, by considering the homotopy exact sequence of
$f_{K,2^V}$.

Corollary \ref{commsub} confirms that
\begin{equation}\label{nicebun}
[G,L]\longrightarrow\mathcal{Z}_G(K,L)\stackrel{p}{\longrightarrow}
\mathcal{Z}_G(K,2^V)
\end{equation}
is a principal $[G,L]$-bundle, classified by $\gamma$. This bundle
encodes a wealth of geometrical information on the pair $(L,K)$. Its
total space measures the failure of $f_{K,L}$ to be a homotopy
equivalence, and its base space is the complement of the coordinate
subspace arrangement $\mathcal{A}_\bF(K)$ by Corollary \ref{hofili}.
Moreover, Theorem \ref{lkequiv} implies that $\gamma$ is also a
$d(K,L)$-equivalence, and so sheds some light on the homotopy type of
$U_\bF(K)$.

Looping \eqref{diago} gives a homotopy commutative diagram of fibrations
\begin{equation}\label{diagt}
\begin{CD}
\varOmega\mathcal{Z}_G(K,L)@>id>>\varOmega\mathcal{Z}_G(K,L)@>>>1\\
@VV\varOmega pV@VVV@VVV\\
\varOmega U_\bF(K)@>i>>\varOmega\djs(K)@>
\;\varOmega f_{K,2^V}\;>>G^m\\
@VV\varOmega\gamma V@VV\varOmega f_{K,L}V@VVidV\\
[G,L]@>>>\coltmg G^L@>a>>G^m
\end{CD}
\end{equation}
in \cat{tmonh}, which offers an alternative perspective on
$\varOmega\djs(K)$.
\begin{lem}\label{lssplit}
The loop space $\varOmega\djs(K)$ splits as $G^m\times\varOmega
U_\bF(K)$ for any simplicial complex $K$; the splitting is not
multiplicative.
\end{lem}
\begin{proof}
The vertex groups $G_j$ embed in $\varOmega\djs(K)$ via homotopy
homomorphisms, whose product $j\colon G^m\rightarrow\varOmega\djs(K)$ is
left inverse to $\varOmega f_{K,2^V}$ (but not a homotopy
homomorphism). The product of the maps $i$ and $j$ is the required
homeomorphism.
\end{proof}

The following examples continue the theme of Examples \ref{exast}
and \ref{arrncomps}. They refer to the second horizontal fibration
of the diagram \eqref{diagt}, which is homotopy equivalent to the third
whenever $K=L$ is flag, by Proposition \ref{hkkcontract}. The second
examples also appeal to James's Theorem \cite{ja:rps}, which
identifies the loop space $\varOmega S^n$ with the free monoid
$F^+(S^{n-1})$ for any $n>1$.
\begin{exas}\label{exasth}
If $K$ is the discrete flag complex $V$, then $\varOmega U_\bF(K)$ is
homotopy equivalent to the commutator subgroup of the free product
$\bigfree_{j=1}^mG_j$. If $K$ is the non-flag complex $\partial(m)$,
then $\varOmega U_\bF(K)$ is homotopy equivalent to $F^+(S^{m-2})$ for
$G=\bZ/2$, and $F^+(S^{2m-2})$ for $G=T$; the map $i$ identifies
the generators of each free monoid with higher Samelson products
(of order $m$) in $\varOmega\djs(K)$.
\end{exas}
Of course, both examples split topologically according to Lemma
\ref{lssplit}. The appearance of higher products in
$\varOmega\djs(\partial(m))$ shows that commutators alone cannot model
$\varOmega\djs(K)$ when $K$ is not flag. More subtle structures are
required, based on higher homotopy commutativity; they are related to
Samelson and Whitehead products, as we explain elsewhere
\cite{para:hhc}.

\section{Homotopy colimits of topological monoids}\label{hocotomo}

We now turn to the loop space $\varOmega\djs(K)$ for a general
simplicial complex $K$, appealing to the theory of homotopy
colimits. Although the resulting models are necessarily more
complicated, they are homotopy equivalent to $\coltmg G^K$ when $K$ is
flag. The constructions depend fundamentally on the categorical ideas
of Section \ref{capr}, and apply to more general spaces than
$\djs(K)$. We therefore work with an arbitrary diagram
$D\colon\cat{a}\rightarrow\cat{tmg}$ for most of the section, and
write $BD\colon\cat{a}\rightarrow\cat{top$_+$}$ for its classifying
diagram. Our applications follow by substituting $G^K$ for $D$.

We implement proposals of earlier authors (as in \cite{vo:hlc}, for
example) by forming the homotopy colimit $\hoc^{\scat{tmg}}D$ in
$\cat{tmg}$, rather than $\cat{top}_+$. This is made possible by the
observation of Section \ref{capr} that the categories $\cat{tmg}$ are
$\cat{t}$-cocomplete, and therefore have sufficient structure for the
creation of internal homotopy colimits. We confirm that $\hoctmg D$ is
a model for the loop space $\varOmega\hoc^+BD$ by proving that $B$
commutes with homotopy colimits in the relevant sense. As usual, we
work in $\cat{tmg}$, but find it convenient to describe certain
details in terms of topological monoids; whenever these monoids are
topological groups, so is the output.

We recall the standard extension of the 2-sided bar construction to the
based setting, with reference to \eqref{tsbar}. We write
$B_\bullet^+(*,\cat{a},D)$ for the diagram
$\cat{b}^{op}\times\del^{op}\rightarrow\cat{top}_+$ given by
\[
(b,(n))\longmapsto
\bigvee_{a_0,a_n}D(b,a_0)\wedge\cat{a}_n(a_0,a_n)_+,
\]
where $D$ is a diagram
$\cat{a$\times$b$^{op}$}\rightarrow\cat{top}_+$. Following Examples
\ref{defhoco}, we define the homotopy $\cat{top}_+$-colimit as
\[
\hoc^+D=B^+(*,\cat{a},D),
\]
and note the equivalent expressions
$B^+(*,\cat{a},A_+)\wedge_{\scat{a}}D\cong
D\wedge_{\scat{a}^{op}}B^+(*,\cat{a},A_+)$.

For $\cat{tmg}$, we proceed by categorical analogy. We replace the
$\cat{top}$-coproduct in \eqref{tsbar} by its counterpart in
$\cat{tmg}$, and the internal cartesian product in $\cat{top}$ by the
tensored structure of $\cat{tmg}$ over $\cat{top}_+$. For any diagram
$D\colon\cat{a}\rightarrow\cat{tmg}$, the simplicial topological
monoid $B_\bullet^{\scat{tmg}}(*,\cat{a},D)$ is therefore given by
\begin{equation}\label{btmg}
(n)\longmapsto\bigfree_{a_0,a_n}D(a_0)\circledast\cat{a}_n(a_0,a_n)_+,
\end{equation}
where $*$ denotes the free product of topological monoids. The face
and degeneracy operators are defined as before, but are now
homomorphisms. When $\cat{a}$ is of the form $\cat{cat}(K)$, the
$n$-simplices \eqref{btmg} may be rewritten as the finite free product
\[
B^{\scat{tmg}}_n(*,\cat{cat}(K),D)=
\bigfree_{\sigma_n\supseteq\cdots\supseteq\sigma_0}D(\sigma_0),
\]
where there is one factor for each $n$-chain of simplices in $K$.
\begin{defn}\label{tmghoco}
The {\it homotopy \cat{tmg}-colimit\/} of $D$ is given by
\[
\hoctmg D=
|B_{\bullet}^{\scat{tmg}}(*,\cat{a},D)|_{\scat{tmg}}
\]
in $\cat{tmg}$, for any diagram $D\colon\cat{a}\rightarrow\cat{tmg}$.
\qed
\end{defn}
So $\hoctmg D$ is an object of $\cat{tmg}$. Following Construction
\eqref{tmontens}, it may be described in terms of generators and
relations as a quotient monoid of the form
\[
\Big(\bigfree_{n\geq 0}\left(B_n(*,\cat{a},D)\circledast\varDelta^n_+
\right)\Big)\;\Big/\;
\Big\langle
\big(d^i_n(b),s\big)=\big(b,\delta^i_n(s)\bigr),\;
\bigl(s^i_n(b),t\bigr)=\bigl(b,\sigma^i_n(t)\bigr)
\Big\rangle\;,
\]
for all $b\in B_n(*,\cat{a},D)$, and all $s\in\varDelta(n-1)$ and
$t\in\varDelta(n+1)$. Here $\delta^i_n$ and $\sigma^i_n$ are the
standard face and degeneracy maps of geometric simplices.
\begin{exa}\label{monmapcyl}
Suppose that $\cat{a}$ is the category
$\,\cdot\rightarrow\cdot\,$, with a single non-identity. Then an
$\cat{a}$-diagram is a homomorphism $M\rightarrow N$ in
$\cat{tmg}$, and $\hoctmg D$ is its $\cat{tmg}$ mapping
cylinder. It may be identified with the $\cat{tmg}$-pushout of the
diagram
\[
M\circledast\varDelta(1)_+\stackrel{j}{\longleftarrow}M
\longrightarrow N,
\]
where $j(m)=(m,0)$ in $M\circledast\varDelta(1)_+$ for all $m\in M$.
\end{exa}

An alternative expression for the simplicial topological monoid
$B_\bullet^{\scat{tmg}}(*,\cat{a},D)$ arises by analogy with the
equivalences \eqref{bars}.
\begin{prop}\label{eqdef}
There is an isomorphism
$D\circledast_{\scat{a}^{op}}B^+_\bullet(*,\cat{a},A_+)\cong
B^{\scat{tmg}}_{\bullet}(*,\cat{a},D)$ of simplicial
topological monoids, for any diagram
$D\colon\cat{a}\rightarrow\cat{tmg}$.
\end{prop}
\begin{proof}
By \eqref{duht}, the functor
$D\circledast_{\scat{a}^{op}}\vble\colon
\fcat{a$^{op}\times$\del$^{op}$}{top$_+$}\rightarrow
\fcat{\del$^{op}$}{tmg}$ is left $\cat{top}_+$-adjoint to
$\cat{tmg}(D,\vble)$, and therefore preserves coproducts. So we may
write
\begin{align*}
D\circledast_{\scat{a}^{op}}B_{\bullet}^+(*,\cat{a},A_+)
&\cong\bigfree_{a,b}D\circledast_{\scat{a}^{op}}
\left(\cat{a}(\vble,a)_+\wedge\cat{a}_\bullet(a,b)_+\right)\\
&\cong\bigfree_{a,b}D(a)\circledast\cat{a}_{\bullet}(a,b)_+
\end{align*}
as required, using the isomorphism
$D\circledast_{\scat{a}^{op}}\cat{a}(\vble,a)\cong D(a)$ of
\eqref{asscot}.
\end{proof}

It is important to establish when the simplicial topological monoids
$B^{\scat{tmg}}_{\bullet}(*,\cat{a},D)$ are proper simplicial spaces, in
the sense of \cite{ma:esg}, because we are interested in the homotopy
type of their realisations. This is achieved in Proposition
\ref{btmgprop}, and leads on to the analogue of the Homotopy Lemma for
$\cat{tmg}$. These are two of the more memorable of the following
sequence of six preliminaries, which precede the proof of our main
result. On several occasions we insist that objects of $\cat{tmg}$ are
well-pointed, and even that they have the homotopy type of a
CW-complex. Such conditions certainly hold for our exponential diagrams,
and do not affect the applications.

We consider families of monoids indexed by the elements $s$ of an
arbitrary set $S$.
\begin{lem}\label{cophes}
Let $f_s\colon M_s\rightarrow N_s$ be a family of homomorphisms of
well-pointed monoids, which are homotopy equivalences; then the
coproduct homomorphism
\[
\bigfree_sf_s\colon\bigfree_sM_s\longrightarrow\bigfree_sN_s
\]
is also a homotopy equivalence.
\end{lem}
\begin{proof}
Let $f\colon M\rightarrow N$ denote the homomorphism in question, and
write $F_kM$ for the subspace of $M$ of elements representable by words
of length $\leq k$. Hence $F_0=\{e\}$, and $F_{k+1}M$ is the pushout
\begin{equation}\label{pusho}
\begin{CD}
\bigvee_KW_K(M)@>\subseteq>>\bigvee_K P_K(M)\\
@VVV@VVV\\
F_kM@>j_k>>F_{k+1}M
\end{CD}
\end{equation}
in $\cat{top}_+$, where $K$ runs through all $(k+1)$-tuples
$(s_1,\dots,s_{k+1})\in S^{k+1}$ such that $s_{i+1}\neq s_i$, and
$W_K(M)\subset P_K(M)$ is the fat wedge subspace of
$M_{s_1}\times\ldots\times M_{s_{k+1}}$. Each $M_s$ is well-pointed, so
$W_K(M)\subset P_K(M)$ is a closed cofibration, and therefore so is
$j_k$. Since $M=\col_k F_k M$ in $\cat{top}_+$, it remains to confirm
that the restriction $f_k\colon F_kM\rightarrow F_kN$ is a homotopy
equivalence for all $k$. We proceed by induction, based on the trivial
case $k=0$.

The map $f$ induces a homotopy equivalence $W_K(M)\rightarrow W_K(N)$
because $M_s$ and $N_s$ are well-pointed, and a further homotopy
equivalence $P_K(M)\rightarrow P_K(N)$ by construction. So the inductive
hypothesis combines with Brown's Gluing Lemma \cite[2.4]{wezizi:hcc} to
complete the proof.
\end{proof}
\begin{lem}\label{subcop}
For any subset $R\subset S$, the inclusion $\ast_rM_r\rightarrow
\ast_sM_s$ is a closed cofibration; in particular, $\ast_sM_s$ is
well-pointed.
\end{lem}
\begin{proof}
Let $B\rightarrow M$ be the inclusion in question, with $F_kM$ as in the
proof of Lemma \ref{cophes} and $F'_kM =B\cup F_kM$. Then $F'_{k+1}M$ is
obtained from $F'_kM$ by attaching spaces $P_K(M)$, where $K$ runs
through all $(s_1,\ldots,s_{k+1})$ in $S^{k+1}\setminus R^{k+1}$ such
that $s_{i+1}\neq s_i$. Thus
$B\rightarrow F'_kM$ is a cofibration for all $k$, implying the result.
\end{proof}

\begin{prop}\label{btmgprop}
Given any small category $\cat{a}$, and any diagram
$D\colon\cat{a}\rightarrow\cat{tmg}$ of well-pointed monoids, the
simplicial space $B^{\scat{tmg}}_{\bullet}(\ast,\cat{a},D)$ is proper,
and its realisation $B^{\scat{tmg}}(\ast,\cat{a},D)$ is well-pointed.
\end{prop}
\begin{proof}
By Lemma \ref{subcop}, each degeneracy map
$B^{\scat{tmg}}_n(\ast,\cat{a},D)\rightarrow
B^{\scat{tmg}}_{n+1}(\ast,\cat{a},D)$ is
a closed cofibration. The first result then follows from Lillig's Union
Theorem \cite{li:utc} for cofibrations. So
$B^{\scat{tmg}}_0(\ast,\cat{a},D)\subset B^{\scat{tmg}}(\ast,\cat{a},D)$
is a closed cofibration and $B^{\scat{tmg}}_0(\ast,\cat{a},D)$ is
well-pointed, yielding the second result.
\end{proof}

As described in Examples \ref{tencotexa}, every simplicial object
$M_\bullet$ in $\cat{tmg}$ has two possible realisations. We now confirm
that they agree, and identify their classifying space.
\begin{lem}\label{tworeals}
The realisations $|M_\bullet|_{\scat{tmg}}$ and $|M_\bullet|$ are
naturally isomorphic objects of $\cat{tmg}$, whose classifying space is
naturally homeomorphic to $|B(M_\bullet)|$.
\end{lem}
\begin{proof}
We apply the techniques of \cite[VII \S3]{elkrmama:rma} and
\cite[\S4]{mcscvo:thh} to the functors $|\vble|_{\scat{tmg}}$ and the
restriction of $|\vble|$ to $\fcat{\del$^{op}$}{tmg}$. Both are left
$\cat{top}_+$-adjoint to
$\Sin\colon\cat{tmg}\rightarrow\fcat{\del$^{op}$}{tmg}$, and so are
naturally equivalent. The homeomorphism $B|M_\bullet|\cong
|B(M_\bullet)|$ arises by considering the bisimplicial object
$(k,n)\mapsto(M_n)^k$ in $\cat{top}_+$, and forming its realisation in
either order.
\end{proof}

We may now establish our promised Homotopy Lemma.
\begin{prop}\label{htylemtmg}
Given diagrams $D_1$, $D_2\colon\cat{a}\rightarrow\cat{tmg}$ of
well-pointed topological monoids, and a map $f\colon
D_1\rightarrow D_2$ such that $f(a)\colon D_1(a)\rightarrow D_2(a)$ is a
homotopy equivalence of underlying spaces for each object $a$ of
$\cat{a}$, the induced map
\[
\hoctmg D_1\longrightarrow\hoctmg D_2
\]
is a homotopy equivalence.
\end{prop}
\begin{proof}
This follows directly from Lemmas \ref{cophes} and \ref{tworeals}, and
Proposition \ref{btmgprop}.
\end{proof}

We need one more technical result concerning homotopy limits of
simplicial objects. We work with diagrams
$X_\bullet\colon\cat{a}\times\del^{op}\rightarrow\cat{top}_+$ of
simplicial spaces, and
$D_\bullet\colon\cat{a}\times\del^{op}\rightarrow\cat{tmg}$ of
simplicial topological monoids.
\begin{prop}\label{twohocs}
With $X_\bullet$ and $D_\bullet$ as above, there are natural
isomorphisms
\[
\hoc^+|X_\bullet|\cong|\hoc^+X_\bullet|\sands
\hoctmg|D_\bullet|\cong|\hoctmg D_\bullet|
\]
in $\cat{top}_+$ and $\cat{tmg}$ respectively.
\end{prop}
\begin{proof}
The isomorphisms arise from realising the bisimplicial objects
\[
(k,n)\longmapsto B^+_k(\ast,\cat{a},X_n)\spandsp (k,n)\longmapsto
B^{\scat{tmg}}_k(\ast,\cat{a},D_n)
\]
in either order. In the case $D_\bullet$, we must also apply the first
statement of Lemma \ref{tworeals}.
\end{proof}

Parts of the proofs above may be rephrased using variants of the
equivalences \eqref{asscot}. They lead to our first general result,
which states that the formation of classifying spaces commutes with
homotopy colimits in an appropriate sense.
\begin{thm}\label{gdishe}
For any diagram $D\colon\cat{a}\rightarrow\cat{tmg}$ of well-pointed
topological monoids with the homotopy types of CW-complexes, the
map
\[
g_D\colon\hoc^+BD\longrightarrow B\hoctmg D
\]
is a homotopy equivalence.
\end{thm}
\begin{proof}
For each object $a$ of $\cat{a}$, let $D_\bullet(a)$ be the singular
simplicial monoid of $D(a)$. The natural map
$|D_\bullet(a)|\rightarrow D(a)$ is a homomorphism of well-pointed
monoids and a homotopy equivalence, so it passes to a homotopy
equivalence $B|D_\bullet(a)|\rightarrow BD(a)$ under the formation
of classifying spaces. By Proposition \ref{htylemtmg} and the
corresponding Homotopy Lemma for $\cat{top}_+$, it therefore
suffices to prove our result for diagrams of realisations of
simplicial monoids.

So let $D_\bullet\colon\cat{a}\times\del^{op}\rightarrow\cat{tmg}$
be a diagram of discrete simplicial monoids. By Lemma \ref{tworeals}
and Proposition \ref{twohocs}, we must show that the canonical map
\[
|\hoc^+BD_\bullet|\longrightarrow|B\hoctmg D_\bullet|
\]
is a homotopy equivalence. Since the simplicial spaces
$\hoc^+BD_\bullet$ and $B\hoctmg D_\bullet$ are proper, this reduces
to proving that
\[
\hoc^+BD_n\longrightarrow B\hoctmg D_n
\]
is a homotopy equivalence in each dimension. But $\hoc^+BD_n$ is the
realisation of the proper simplicial space
$B_\bullet^+(\ast,\cat{a},BD_n)$, and $B\hoctmg D_n$ is naturally
homeomorphic to the proper simplicial space
$B(B_\bullet^{\scat{tmg}}(\ast,\cat{a},D_n))$ by Lemma
\ref{tworeals}. Moreover, for each $k\geq 0$ the natural map
$B_k^+(\ast,\cat{a},BD_n)\rightarrow
B(B_k^{\scat{tmg}}(\ast,\cat{a},D_n)$ coincides with the map
\begin{equation}\label{fiedo}
\bigvee_{a_0\rightarrow\dots\rightarrow a_k} BD_n(a_0)
\longrightarrow B\big(\bigfree_{a_0\rightarrow\dots\rightarrow a_k}
D_n(a_0)\big)
\end{equation}
induced by the inclusion of each $D_n(a_0)$ into the free product.
Since \eqref{fiedo} is a homotopy equivalence by a theorem of
Fiedorowicz \cite[4.1]{fi:cst}, the proof is complete.
\end{proof}

Various steps in the proof of Theorem \ref{gdishe} may be adapted to
verify the following, which answers a natural question about
tensored monoids.
\begin{prop}\label{classtens} For any well-pointed topological
monoid $M$ and based space $Y$, the natural map
\[
BM\wedge Y\longrightarrow B(M\otimes Y)
\]
is a homotopy equivalence when $M$ and $Y$ have the homotopy type of
CW-complexes.
\end{prop}
\begin{proof}
As in Theorem \ref{gdishe}, we need only work with the realisations
$|M_\bullet|$ and $|Y_\bullet|$ of the total singular complexes.
Since $B|M_\bullet|\wedge|Y_\bullet|\rightarrow
B(|M_\bullet|\otimes|Y_\bullet|)$ is the realisation of the natural
map $BM_n\wedge Y_n\rightarrow B(M_n\otimes Y_n)$, it suffices to
assume that $Y$ is discrete; in this case,
\[
BM\wedge Y\longrightarrow B\big(\bigfree_y M_y\big)
\]
is a homotopy equivalence by the same result of Fiedorowicz
\cite{fi:cst}.
\end{proof}

We apply Theorem \ref{gdishe} to construct our general model for
$\varOmega\djs(K)$, but require a commutative diagram to clarify its
relationship with the special case $h_K$ of Proposition \ref{dequiv}.
We deal with $\cat{a}^{op}\times\del^{op}$-diagrams $X_\bullet$ in
$\cat{top}_+$, and certain of their morphisms. These include
$\theta\colon X_\bullet\rightarrow
\cat{top}_+(BD,B(D\circledast_{\scat{a}^{op}}X_\bullet))$, defined for
any $X_\bullet$ by $\theta(x)=B(d\mapsto d\circledast x)$, and the
projection $\pi\colon B^+_\bullet\rightarrow(*_+)_\bullet$, where
$B^+_\bullet$ and $(*_+)_\bullet$ denote $B^+_\bullet(*,\cat{a},A)$ and
the trivial diagram respectively.  Under the homeomorphism
\[
\fcat{\del$^{op}$}{top$_+$}\left(BD\wedge_{\scat{a}^{op}}X_\bullet,
B(D\circledast_{\scat{a}^{op}}X_\bullet)\right)\cong
\fcat{a$^{op}\times$\del$^{op}$}{top$_+$}\big(X_\bullet,
\cat{top$_+$}(BD,B(D\circledast_{\scat{a}^{op}}X_\bullet))\big)
\]
of \eqref{duht}, $\theta$ corresponds to a map $\phi\colon
BD\wedge_{\scat{a}^{op}}X_\bullet\rightarrow
B(D\circledast_{\scat{a}^{op}}X_\bullet)$ of simplicial spaces.
\begin{prop}\label{bhoc}
For any diagram $D\colon\cat{a}\rightarrow\cat{tmg}$, there is a
commutative square
\[
\begin{CD}
\hoc^+ BD@>g_D>>B\hoctmg D\\
@VVp^{\smash{+}}V@VVBp^{\smash{\scat{tmg}}}V\\
\col^+BD@>f_D>>B\coltmg D
\end{CD},
\]
where $p^+$ and $p^{\scat{tmg}}$ are the natural projections.
\end{prop}
\begin{proof}
By construction, the diagram
\[
\begin{CD}
B^+_\bullet@>\theta>>
\cat{top$_+$}\left(BD,B(D\circledast_{\scat{a}^{op}}B^+_\bullet)\right)\\
@VV\pi V@VVB(1\circledast\pi)\cdot V\\
(*_+)_\bullet@>\theta>>
\cat{top$_+$}\left(BD,B(D\circledast_{\scat{a}^{op}}(*_+)_\bullet)\right)
\end{CD}
\]
is commutative in $\fcat{a$^{op}\times$\del$^{op}$}{top$_+$}$, and has
adjoint
\begin{equation}\label{helpdi}
\begin{CD}
BD\wedge_{\scat{a}^{op}}B^+_\bullet@>\phi>>
B(D\circledast_{\scat{a}^{op}}B^+_\bullet)\\
@VV1\wedge\pi V@VVB(1\circledast\pi) V\\
BD\wedge_{\scat{a}^{op}}(*_+)_\bullet@>\phi>>
B(D\circledast_{\scat{a}^{op}}(*_+)_\bullet)
\end{CD}
\end{equation}
in $\fcat{\del$^{op}$}{top$_+$}$. By Proposition \ref{eqdef}, the upper
$\phi$ is the map $B_{\bullet}^+(*,\cat{a},BD)\rightarrow
B\left(B_{\bullet}^{\scat{tmg}}(*,\cat{a},D)\right)$ obtained by
applying the relevant map \eqref{fiedo} in each dimension. By Examples
\ref{tencotexa}, the lower $\phi$ is given by the canonical map
$f_D\colon\col^+BD\rightarrow B\coltmg D$ in each dimension. Since
realisation commutes with $B$, the topological realisation of
\eqref{helpdi} is the diagram we seek; for Lemma \ref{tworeals}
identifies the upper right-hand space with $B\hoctmg D$, and Examples
\ref{defhoco} confirms that the vertical maps are the natural
projections.
\end{proof}
\begin{thm}\label{genldj}
There is a homotopy commutative square
\[
\begin{CD}
\varOmega\hoc^+(BG)^K@>\overline{h}_K>>\hoctmg G^K\\
@VV{\smash{\varOmega}}p_KV@VVq_KV\\
\varOmega\djs(K)@>h_K>>\coltmg G^K
\end{CD}
\]
of homotopy homomorphisms, where $p_K$ and $\overline{h}_K$ are homotopy
equivalences for any simplicial complex $K$.
\end{thm}
\begin{proof}
We apply Proposition \ref{bhoc} with $D=G^K$, and loop the corresponding
square; the projection $p_K\colon\hoc^+(BG)^K\rightarrow\djs(K)$ is
a homotopy equivalence, as explained in Corollary \ref{hofibkl}.
The result follows by composing the horizontal maps with the
canonical homotopy homomorphism $\varOmega BH\rightarrow H$, where
$H=\hoctmg G^K$ and $\coltmg G^K$ respectively.
\end{proof}

It is an interesting challenge to describe good geometrical models for
homotopy homomorphisms which are inverse to $h_K$ and $\overline{h}_K$.


\begin{thebibliography}{10}

\bibitem{bawe:ttt}
Michael Barr and Charles Wells.
\newblock {\em Toposes, Triples and Theories}.
\newblock Number 278 in Grundlehren der mathematischen Wissenschaften.
  Springer-Verlag, 1985.

\bibitem{bj:sa}
Anders Bj{\"o}rner.
\newblock Subspace arrangements.
\newblock In A~Joseph, F~Mignot, F~Murat, and B~Prum, editors, {\em Proceedings
  of the First European Congress of Mathematics (Paris 1992)}, volume 119 of
  {\em Progress in Mathematics}, pages 321--370. Birkh{\"a}user, 1995.

\bibitem{bo:hca}
Francis Borceux.
\newblock {\em Handbook of Categorical Algebra 2}, volume~51 of {\em
  Encyclopedia of Mathematics and its Applications}.
\newblock Cambridge University Press, 1994.

\bibitem{boka:hlc}
A~K Bousfield and Daniel~M Kan.
\newblock {\em Homotopy Limits, Completions and Localizations}, volume 304 of
  {\em Lecture notes in Mathematics}.
\newblock Springer Verlag, 1972.

\bibitem{br:shc}
Michael Brinkmeier.
\newblock Strongly homotopy-commutative monoids revisited.
\newblock {\em Documenta Mathematica}, 5:613--624, 2000.

\bibitem{brha:tgo}
Ronald Brown and J~P~Lewis Hardy.
\newblock Topological groupoids 1: universal constructions.
\newblock {\em Mathematische Nachriften}, 71:273--286, 1976.

\bibitem{bupa:tact}
Victor~M Buchstaber and Taras~E Panov.
\newblock Torus actions, combinatorial topology, and homological algebra.
\newblock {\em Russian Mathematical Surveys}, 55:825--921, 2000.

\bibitem{cafo:pcc}
Pierre Cartier and Dominique Foata.
\newblock {\em Probl{\`e}mes combinatoires de commutation et
  r{\'e}arrangements}, volume~85 of {\em Lecture Notes in Mathematics}.
\newblock Springer Verlag, 1969.

\bibitem{chda:fka}
Ruth Charney and Michael Davis.
\newblock Finite {$K(\pi,1)$} for {A}rtin groups.
\newblock In {\em Prospects in Topology (Princeton NJ, 1994)}, volume 138 of
  {\em Annals of Mathematics Studies}, pages 110--124. Princeton University
  Press, 1995.

\bibitem{ch:ecg}
Ian~M Chiswell.
\newblock The {E}uler characteristic of graph products and of {C}oxeter groups.
\newblock In William~J Harvey and Colin Maclachlan, editors, {\em Discrete
  Groups and Geometry (Birmingham, 1991)}, volume 173 of {\em London
  Mathematical Society Lecture Note Series}, pages 36--46. Cambridge University
  Press, 1992.

\bibitem{da:ggr}
Michael~W Davis.
\newblock Groups generated by reflections and aspherical manifolds not covered
  by {E}uclidean space.
\newblock {\em Annals of Mathematics}, 117:293--324, 1983.

\bibitem{daja:cpc}
Michael~W Davis and Tadeusz Januszkiewicz.
\newblock Convex polytopes, {C}oxeter orbifolds and torus actions.
\newblock {\em Duke Mathematical Journal}, 62:417--451, 1991.

\bibitem{dola:hsp}
Albrecht Dold and Richard Lashof.
\newblock Principal quasifibrations and fibre homotopy equivalence of bundles.
\newblock {\em Illinois Journal of Mathematics}, 3:285--305, 1959.

\bibitem{elkrmama:rma}
Anthony~D Elmendorf, Igor Kriz, Michael~P Mandell, and J~Peter May.
\newblock {\em Rings, Modules, and Algebras in Stable Homotopy Theory},
  volume~47 of {\em Mathematical Surveys and Monographs}.
\newblock American Mathematical Society, 1997.

\bibitem{fi:cst}
Zbigniew Fiedorowicz.
\newblock Classifying spaces of topological monoids and categories.
\newblock {\em Journal of the American Mathematical Society}, 106:301--350,
  1984.

\bibitem{gule:hlr}
T~Gulliksen and G~Levin.
\newblock {\em Homology of local rings}, volume~20 of {\em Queen's papers in
  Pure and applied Mathematics}.
\newblock Queen's University, Kingston, Ontario, 1969.

\bibitem{herewe:cli}
J{\o}rgen Herzog, Victor Reiner, and Volkmar Welker.
\newblock Componentwise linear ideals and {G}olod rings.
\newblock {\em Michigan Mathematical Journal}, 46:211--223, 1999.

\bibitem{hovo:mts}
Jens Hollender and Rainer~M Vogt.
\newblock Modules of topological spaces, applications to homotopy limits and
  {$E_\infty$} structures.
\newblock {\em Archiv der Mathematik}, 59:115--129, 1992.

\bibitem{ja:rps}
Ioan~M James.
\newblock Reduced product spaces.
\newblock {\em Annals of Mathematics}, 62:259--280, 1955.

\bibitem{joro:cbc}
S~A Joni and G-C Rota.
\newblock Coalgebras and bialgebras in combinatorics.
\newblock {\em Studies in Applied Mathematics}, 61:93--139, 1979.

\bibitem{ke:bce}
G~Maxwell Kelly.
\newblock {\em Basic Concepts of Enriched Category Theory}, volume~64 of {\em
  London Mathematical Society Lecture Note Series}.
\newblock Cambridge University Press, 1982.

\bibitem{kiro:hca}
Ki~Hang Kim and Fred~W Roush.
\newblock Homology of certain algebras defined by graphs.
\newblock {\em Journal of Pure and Applied Algebra}, 17:179--186, 1980.

\bibitem{li:utc}
Joachim Lillig.
\newblock A union theorem for cofibrations.
\newblock {\em Archiv der Mathematik}, 24:410--415, 1973.

\bibitem{ma:soa}
J~Peter May.
\newblock {\em Simplicial Objects in Algebraic Topology}, volume~11 of {\em Van
  Nostrand Mathematical Studies}.
\newblock Van Nostrand Reinhold, 1967.

\bibitem{ma:esg}
J~Peter May.
\newblock {$E_\infty$}-spaces, group completions and permutative categories.
\newblock In Graeme Segal, editor, {\em New Developments in Topology},
  volume~11 of {\em London Mathematical Society Lecture Notes Series}, pages
  153--231. Cambridge University Press, 1974.

\bibitem{mcscvo:thh}
James McClure, Roland Schw{\"a}nzl, and Rainer Vogt.
\newblock {\it THH$(R)\cong R\otimes S^1$} for {$E_\infty$} ring spectra.
\newblock {\em Journal of Pure and Applied Algebra}, 121:137--159, 1997.

\bibitem{para:hhc}
Taras Panov and Nigel Ray.
\newblock The homology and homotopy theory of certain loop spaces.
\newblock In preparation, University of Manchester, 2001.

\bibitem{perewe:crd}
Irena Peeva, Victor Reiner, and Volkmar Welker.
\newblock Cohomology of real diagonal subspace arrangments via resolutions.
\newblock {\em Compositio Mathematica}, 117:99--115, 1999.

\bibitem{rasc:cmc}
Nigel Ray and William Schmitt.
\newblock Combinatorial models for coalgebraic structures.
\newblock {\em Advances in Mathematics}, 138:211--262, 1998.

\bibitem{scvo:cae}
Roland Schw{\"a}nzl and Rainer~M Vogt.
\newblock The categories of {$A_\infty$} and {$E_\infty$} monoids and ring
  spectra as closed simplicial and topological model categories.
\newblock {\em Archiv der Mathematik}, 56:405--411, 1991.

\bibitem{wa:tmg}
Peter Scott and Terry Wall.
\newblock Topological methods in group theory.
\newblock In C~T~C Wall, editor, {\em Homological Group Theory (Proceedings of
  the Durham Symposium, 1977)}, volume~36 of {\em London Mathematical Society
  Lecture Note Series}, pages 137--203. Cambridge University Press, 1979.

\bibitem{se:css}
Graeme Segal.
\newblock Classifying spaces and spectral sequences.
\newblock {\em Institut des Hautes {\'E}tudes Scientifiques, Publications
  Mathematiques}, 34:105--112, 1968.

\bibitem{st:cca}
Richard~P Stanley.
\newblock {\em Combinatorics and Commutative Algebra, 2nd edition}, volume~41
  of {\em Progress in Mathematics}.
\newblock Birkh{\"a}user, Boston, 1996.

\bibitem{su:hcg}
Masahiro Sugawara.
\newblock On the homotopy-commutativity of groups and loop spaces.
\newblock {\em Memoirs of the College of Science, University of Kyoto: Series
  A, Mathematics}, 33:257--269, 1960/61.

\bibitem{vo:cct}
Rainer~M Vogt.
\newblock Convenient categories of topological spaces for homotopy theory.
\newblock {\em Archiv der Mathematik}, 22:545--555, 1971.

\bibitem{vo:hlc}
Rainer~M Vogt.
\newblock Homotopy limits and colimits.
\newblock {\em Mathematische Zeitschrift}, 134:11--52, 1973.

\bibitem{wezizi:hcc}
Volkmar Welker, G{\"u}nter~M Ziegler, and Rade~T {\v Z}ivaljevi{\'c}.
\newblock Homotopy colimits - comparison lemmas for combinatorial applications.
\newblock {\em Journal f{\"u}r die reine und angewandte Mathematik},
  509:117--149, 1999.

\bibitem{zizi:hts}
G{\"u}nter~M Ziegler and Rade~T {\v Z}ivaljevi{\'c}.
\newblock Homotopy types of subspace arrangements via diagrams of spaces.
\newblock {\em Mathematische Annalen}, 295:527--548, 1993.

\end{thebibliography}

\end{document}